\documentclass[11pt,reqno]{amsart}

\usepackage{amsmath}
\usepackage{amssymb}
\usepackage{amsthm}

\newtheorem{theorem}{Theorem}[section]
\newtheorem{prop}[theorem]{Proposition}
\newtheorem{lem}[theorem]{Lemma}
\newtheorem*{cor}{Corollary}
\theoremstyle{definition}
\newtheorem{defn}[theorem]{Definition}
\theoremstyle{remark}
\newtheorem*{rem}{Remark}
\numberwithin{equation}{section}

\begin{document}

\title[Quantum double construction]
{Quantum double construction in the $C^*$-algebra setting of certain
Heisenberg-type quantum groups}

\author{Byung-Jay Kahng}
\date{}
\address{Department of Mathematics and Statistics\\ University of Nevada\\
Reno, NV 89557}
\email{bjkahng@unr.edu}

\begin{abstract}
In this paper, we carry out the {\em quantum double construction\/} of
the specific quantum groups we constructed earlier, namely, the ``quantum
Heisenberg group algebra'' $(A,\Delta)$ and its dual, the ``quantum Heisenberg
group'' $(\hat{A},\hat{\Delta})$.  Our approach will be by constructing
a suitable multiplicative unitary operator. In this way, we are able to
retain the $C^*$-algebra framework, and thus able to carry out our
construction within the category of locally compact quantum groups.
This construction is a kind of a generalized crossed product.

 To establish that the quantum double we obtain is indeed a locally compact
quantum group, we will also discuss the dual of the quantum double and the
Haar weights on both of these double objects.  Towards the end, we also
include a construction of a (quasitriangular) ``quantum universal $R$-matrix''.
\end{abstract}
\maketitle

{\sc Introduction.}
The quantum double construction, which was originally introduced by
Drinfeld in the mid-80's for (finite-dimensional) Hopf algebras \cite{Dr},
is among the most celebrated methods of constructing non-commutative and
non-cocommutative Hopf algebras.  Even in the case of an ordinary group,
equivalently for the algebra of (continuous) functions $C(G)$, the quantum
double construction leads to an interesting crossed product algebra $C(G)
\rtimes_{\alpha}G$, where $\alpha$ is the conjugation \cite{Mj}, \cite{KrM1}.

 We wish to carry out a similar construction in the framework of ($C^*$-algebraic)
locally compact quantum groups.  This is not totally a new endeavor: As early
as in \cite{PW}, Podles and Woronowicz has constructed their example of a quantum
Lorentz group, by considering the quantum double of the compact quantum group
$SU_{\mu}(2)$; Baaj and Skandalis \cite{BS} have a version in the context of
the multiplicative unitary operators; And more recently, Yamanouchi \cite{Ya}
has made this more systematic while Baaj and Vaes \cite{BjVa} considers a
more generalized framework of double crossed products.

 On the other hand, as is the case for a lot of going-ons in the study of
locally compact quantum groups (especially for the non-compact ones), there
have been only a handful of work done on actual examples.  In this paper,
we will obtain the quantum double object of the ``quantum Heisenberg group
algebra'' $(A,\Delta)$ and the ``quantum Heisenberg group'' $(\hat{A},
\hat{\Delta})$, which are the specific non-compact quantum groups we
constructed earlier (See \cite{BJKp2}, \cite{BJKppha}, \cite{BJKqhg}.).
The quantum double will be also a valid locally compact quantum group.

 In addition to finding a new example of a quantum group and enriching
the duality picture between $(A,\Delta)$ and $(\hat{A},\hat{\Delta})$,
there are other merits of studying the quantum double.  An interesting
point is that the quantum double is not just an algebraic object, but
also a nice non-commutative geometric object of study: Note that
since the quantum double is obtained as a generalized crossed product,
it can be considered as a kind of a ``quantized space'', while being a
quantum group means it is also ``group-like''.  It will be an interesting
future research project to further explore how these two different flavors
arise together in our example.

 The goal of the present paper is to give an actual construction of
the quantum double, give a concrete realization as an operator algebra
on a specific Hilbert space, establish it as a $C^*$-algebraic, locally
compact quantum group in the sense of Kustermans and Vaes \cite{KuVa}.
Also constructed are the dual object of the quantum double (also a
locally compact quantum group), and the ``quantum universal $R$-matrix''
type operator for the quantum double.

 Our construction method and techniques are strongly motivated by and
are based on the fundamental paper by Baaj and Skandalis \cite{BS}.
On the other hand, we note that the presentation given in \cite{BS},
as well as the ones in \cite{Ya}, \cite{BjVa}, are somewhat less suitable
for developing a rich connection with the Poisson--Lie group theory.
In our presentation, by explicitly working with a dense subalgebra
of functions contained in a $C^*$-algebra, we make it much easier to 
establish a link between the quantum ($C^*$-algebra) setting and the
classical (Poisson--Lie group) level.  In our future work, we are
planning to further encourage this link and explore the properties
of the quantum double via the tools like Lie bialgebras and dressing
actions.

 Here is how the paper is organized: In Section 1, we briefly summarize
the quantum double construction in the (finite-dimensional) Hopf algebra
setting.  We generally follow Majid \cite{Mj}.  We will use this section
as a guide for our main construction in the $C^*$-algebra framework.

 In Section 2, we describe the specific quantum groups $(A,\Delta)$ and
$(\hat{A},\hat{\Delta})$, which we constructed in our previous papers.
Since the multiplicative unitary operators will play a central role in
the later sections, we chose to give characterizations of $(A,\Delta)$
and $(\hat{A},\hat{\Delta})$ as subalgebras in ${\mathcal B}({\mathcal H})$,
via a multiplicative unitary operator $U_A$.

 Our main construction of the quantum double $D(A)=(A_D,\Delta_D)$ is
carried out in Section~3.  The definition is given in terms of multiplicative
unitary operators, but we provide justification that it is compatible with
the definition in the purely algebraic setting.  Reflecting the fact that
the quantum double construction is closed in the category of ``Kac algebras'',
we note that the antipodal map $S_D$ for our example satisfies ${S_D}^2
\cong\operatorname{Id}$.  In Section 4, we look at the dual of the quantum
double $(\widehat{A_D},\widehat{\Delta_D})$.  Here, $\widehat{A_D}$ is
isomorphic as a $C^*$-algebra to $A\otimes\hat{A}^{\operatorname{op}}$,
but its coalgebra structure is twisted.

 In Section 5, discussion is given on Haar weights for both of the dual
objects $(\widehat{A_D},\widehat{\Delta_D})$ and $(A_D,\Delta_D)$.  We
see that $(\widehat{A_D},\widehat{\Delta_D})$ is unimodular, while $(A_D,
\Delta_D)$ is not.  The existence of the legitimate Haar weights assures
us that both are ($C^*$-algebraic) locally compact quantum groups.

 In Section 6, we find an operator ${\mathcal R}$ contained in the multiplier
algebra $M(A_D\otimes A_D)$, which can be considered as a ``quantum universal
$R$-matrix''.  We only give its construction here.  Its possible applications
to representation theory and connection with the Poisson structure at the
classical limit level will be postponed to a future occasion.

{\sc Terminology.} Let ${\mathcal H}$ be a Hilbert space.  A unitary
operator $V\in{\mathcal B}({\mathcal H}\otimes{\mathcal H})$ is said
to be {\em multiplicative\/}, if it satisfies the following ``pentagon
equation'':
$$
V_{12}V_{13}V_{23}=V_{23}V_{12}\bigl(\in{\mathcal B}({\mathcal H}
\otimes{\mathcal H}\otimes{\mathcal H})\bigr).
$$
Here, the notation $V_{13}$ indicates that the operator $V$ acts
only on the first and third copies of ${\mathcal H}$, while letting
the second copy unchanged.  Similar comments hold for the others.
For a systematic discussion on multiplicative unitary operators,
see the paper by Baaj and Skandalis \cite{BS}.

\section{Quantum double in the purely algebraic framework}

 In this section, we will work only with finite-dimensional Hopf algebras.
The goal here is to collect some useful results from the purely algebraic
setting, which will guide us in our main construction at the level of
($C^*$-algebraic) locally compact quantum groups.  Most of the results
below are standard ones.  See, for instance, \cite{Dr}, \cite{Mo}, \cite{CP},
\cite{Mj}.

 Given a (finite-dimensional) Hopf algebra $B$, the {\em quantum double\/}
$D(B)$ is a certain ``double crossed product'' algebra, $D(B)={B^*}^
{\operatorname{op}}\Join B$, where ${B^*}^{\operatorname{op}}$ is same as
the dual Hopf algebra $B^*$ but equipped with the opposite multiplication.
The Hopf algebras $B$ and ${B^*}^{\operatorname{op}}$ mutually act by
(generalized) coadjoint actions.  A more precise description is given
below.

\begin{defn}\label{qdouble}
\begin{enumerate}
\item The (left) coadjoint action of $B$ on ${B^*}^{\operatorname{op}}$
is defined by
$$
f\triangleright\phi=\operatorname{Ad}^*_f(\phi)=\sum \phi_{(2)}\bigl
\langle f,(S\phi_{(1)})\phi_{(3)}\bigr\rangle.
$$
Similarly, we can define the coadjoint action $\triangleleft$ of $B^*$
on $B$, which we may view as a right action of ${B^*}^{\operatorname{op}}$.
That is,
$$
f\triangleleft\phi=\operatorname{Ad}^*_{\phi}(f)=\sum f_{(2)}\bigl
\langle(Sf_{(1)})f_{(3)},\phi\bigr\rangle.
$$
\item The ``quantum double'' $D(B)={B^*}^{\operatorname{op}}\Join B$ is
equipped with the multiplication:
$$
(\phi\otimes f)\times(\psi\otimes g):=\sum \phi\cdot_{\operatorname{op}}
(f_{(1)}\triangleright\psi_{(1)})\otimes(f_{(2)}\triangleleft\psi_{(2)})g,
$$
and the tensor product comultiplication:
$$
\Delta_{D}(\phi\otimes f):=\sum \phi_{(1)}\otimes f_{(1)}\otimes
\phi_{(2)}\otimes f_{(2)}.
$$
\end{enumerate}
\end{defn}

 In the above, we are using the standard Sweedler notation (See \cite{Mo}.).
That is, we write $\Delta f$ as $\Delta f=\sum f_{(1)}\otimes f_{(2)}$,
and by the coassociativity we have: $(\Delta\otimes\operatorname{id})
(\Delta f)=(\operatorname{id}\otimes\Delta)(\Delta f)=\sum f_{(1)}\otimes 
f_{(2)}\otimes f_{(3)}$.  Since we are considering finite-dimensional
algebras, $\otimes$ denotes the algebraic tensor product.  Meanwhile,
$S$ is the antipode (co-inverse) and $\langle\ ,\ \rangle$ is the
dual pairing.

 The verification of $\triangleright$ and $\triangleleft$ being actions
are not difficult, using the coassociativity and the property of the
antipode map.  As the name suggests, they are generalizations of the
coadjoint actions of groups (similar to taking conjugates).  The
actions make ${B^*}^{\operatorname{op}}$ as a ``$B$-module algebra'',
and $B$ as a ${B^*}^{\operatorname{op}}$-module algebra.  Moreover,
$({B^*}^{\operatorname{op}},B)$ forms a ``matched pair'' of Hopf algebras
(in the sense of Majid), from which the above definition of $D(B)$ arises.
See \cite{Mj1}, \cite{Mj} for details.

 We point out here that above definition of $D(B)$ is different from
Drinfeld's original form \cite{Dr}, containing $B$ and the ``co-opposite
dual'' ${B^*}^{\operatorname{cop}}$.  Ours is actually the one proposed
by Majid (See Theorem 7.1.1 of \cite{Mj}.).  But by means of the antipode
of $B^*$ (considered as a linear map from ${B^*}^{\operatorname{op}}$
to ${B^*}^{\operatorname{cop}}$), we can show without trouble that
the two versions are equivalent.  There are also several other (equally
valid) versions.  Throughout this paper, mainly due to reasons related
with possible future applications, our preferred version of the quantum
double $D(B)$ will be as in Definition \ref{qdouble}.

 It is well known that the quantum double construction leads to a
``quasitriangular'' Hopf algebra.  In case of $D(B)$ as defined
here, its quasitriangular structure is given by $R=\sum_j\psi^j\otimes
f_j$, where $\{f_j\}$ is a basis for $B$ and $\{\psi^j\}$ a dual basis.

 Before we wrap up, let us briefly mention a special case, which will be
a motivating model for us.  Consider an ordinary finite group $G$, and
let $B=\mathbb{C}G$ be the group algebra of $G$ and $B^*=C(G)$ be the
algebra of functions on $G$ (note that $B^*={B^*}^{\operatorname{op}}$
for being commutative), with their natural Hopf algebra structures.  Then
$D(B)$ becomes the crossed product algebra $C(G)\rtimes G$, given by the
conjugate action.  In the case of a locally compact group (not necessarily
finite), this example was studied by Koornwinder and Muller in \cite{KrM1},
\cite{KrM2}.

 Even if this is a rather simple situation, it has an interesting
interpretation as an algebra of quantum observables of a quantum system
(in which a particle is constrained to move on conjugacy classes in $G$).
See Example~6.1.8 of \cite{Mj}.  Other applications can be found in
\cite{DPR}, \cite{BLC}, where the quantum double is used as a generalized
symmetry object.

\section{The quantum Heisenberg group algebra $(A,\Delta)$
and the quantum Heisenberg group $(\hat{A},\hat{\Delta})$}

 We now turn our attention to the $C^*$-algebra setting.  Specifically,
let us consider the non-compact quantum groups $(A,\Delta)$ and $(\hat{A},
\hat{\Delta})$, which were constructed in \cite{BJKp2}, \cite{BJKppha},
\cite{BJKqhg}.  They are mutually dual locally compact quantum groups
(in the sense of Kustermans and Vaes \cite{KuVa}).

 As we pointed out in our previous papers, $(A,\Delta)$ can be regarded
as a ``quantum Heisenberg group algebra'' (i.\,e. ``quantized $C^*(H)$''),
while $(\hat{A},\hat{\Delta})$ may be viewed as a ``quantum Heisenberg
group'' (i.\,e. ``quantized $C_0(H)$'').  Originally, they were obtained
by deformation quantization of the (mutually dual) pair of Poisson--Lie
groups $(G,H)$, where $H$ is the Heisenberg Lie group and $G$ is its dual
Poisson--Lie group carrying a certain non-linear Poisson bracket (See
\cite{BJKp2}, for the description of the non-linear Poisson structure on
$G$ and the construction of $(A,\Delta)$ as its deformation quantization.).  

 However, for the purpose of this article, we will de-emphasize the
deformation process and the role of Poisson geometry.  Instead, our
descriptions of $(A,\Delta)$ and $(\hat{A},\hat{\Delta})$ will be given
in terms of a multiplicative unitary operator $U_A$.  We will postpone
to a separate occasion the discussion of the relationships between the
Poisson-Lie groups $G$, $H$, $D=G\Join H$ (described in \cite{BJKhj})
and the quantum groups $A$, $\hat{A}$, $D(A)$ (to be constructed below).

 Both $C^*$-algebras $A$ and $\hat{A}$ are realized as operator algebras
contained in ${\mathcal B}({\mathcal H})$.  The Hilbert space ${\mathcal H}$
is defined by $L^2(H/Z\times H^*/{Z^{\perp}})$, where $H$ is the
$(2n+1)$-dimensional Heisenberg Lie group (considered naturally as
a vector space), $Z$ is the center of $H$ (which is a subspace of $H$),
while $H^*$ is the dual vector space of $H$, and $Z^{\perp}\subseteq H^*$
is the orthogonal complement of the subspace $Z$.  We think of ${\mathcal H}$
as the space of $L^2$-functions in the $(x,y,r)$ variables, where $(x,y)\in 
H/Z(\cong\mathbb{R}^{2n})$ and $r\in H^*/{Z^{\perp}}(\cong\mathbb{R})$.
By partial Fourier transform in the $r$ variable, ${\mathcal H}$ is isomorphic
to the Hilbert space $L^2(H)$.

 Consider now the unitary operator $U_A\in{\mathcal B}({\mathcal H}\otimes
{\mathcal H})$, as defined in Proposition 3.1 of \cite{BJKp2} (See also
Proposition 2.2 of \cite{BJKppha}):
\begin{align}
U_A\xi(x,y,r;x',y',r')&=(e^{-\lambda r'})^n\bar{e}\bigl[\eta_{\lambda}(r')
\beta(e^{-\lambda r'}x,y'-e^{-\lambda r'}y)\bigr]  \notag \\
&\quad \xi(e^{-\lambda r'}x,e^{-\lambda r'}y,r+r';x'-e^{-\lambda r'}x,
y'-e^{-\lambda r'}y,r').  \notag
\end{align}
Note here that we are using a fairly standard notation of $e(t)=e^{2\pi it}$,
so $\bar{e}(t)=e^{-2\pi it}$.  And $\beta(\ ,\ )$ is the usual inner product.
On the other hand, we need some explanation about the (fixed) constant $\lambda
\in\mathbb{R}$.  It is the constant that determines the aforementioned
non-linear Poisson structure when $\lambda\ne0$ (See \cite{BJKp2}.).  The
expression $\eta_{\lambda}(r)$ is defined such that $\eta_{\lambda}(r)
=\frac{e^{2\lambda r}-1}{2\lambda}$, which reflects the non-linear flavor
(When $\lambda=0$, we take $\eta_{\lambda=0}(r)=r$ which is linear.).

 The unitary operator $U_A$ is multiplicative (satisfying the pentagon equation)
and is regular.  Therefore, by following Baaj and Skandalis \cite{BS}, we can
define a pair of $C^*$-bialgebras $(A,\Delta)$ and $(\hat{A},\hat{\Delta})$.
First, we have:
$$
A=\overline{\bigl\{(\omega\otimes\operatorname{id})(U_A):\omega\in{\mathcal B}
({\mathcal H})_*\bigr\}}^{\|\ \|},
$$
where $L(\omega)=(\omega\otimes\operatorname{id})(U_A)$ are the ``left slices''
of $U_A$ by the linear forms $\omega\in{\mathcal B}({\mathcal H})_*$.  

 For an alternative characterization of $A$, consider ${\mathcal A}$, which is
the space of Schwartz functions in the $(x,y,r)$ variables having compact support
in $r$.  There is the following ``regular representation'' $L$ of ${\mathcal A}$,
on ${\mathcal B}({\mathcal H})$:
$$
(L_f\xi)(x,y,r):=\int f(\tilde{x},\tilde{y},r)\xi(x-\tilde{x},y-\tilde{y},r)
\bar{e}\bigl[\eta_{\lambda}(r)\beta(\tilde{x},y-\tilde{y})\bigr]\,d\tilde{x}
d\tilde{y}.
$$
We have shown in \cite{BJKppha} that $A=\overline{L({\mathcal A})}^{\|\ \|}$.
This means that ${\mathcal A}$ is a (norm dense) ${}^*$-subalgebra of $A$,
and we can regard the functions $f\in{\mathcal A}$ same as the operators $L_f$.
More specifically, the multiplication and the involution on ${\mathcal A}$
take the following form (given by $L_{f\times g}=L_f L_g$ and $L_{f^*}=(L_f)^*$):
\begin{align}\label{(Aops)}
(f\times_A g)(x,y,r)&=\int f(\tilde{x},\tilde{y},r)g(x-\tilde{x},y-\tilde{y},r)
\bar{e}\bigl[\eta_{\lambda}(r)\beta(\tilde{x},y-\tilde{y})\bigr]\,d\tilde{x}
d\tilde{y}.  \notag \\
f^*(x,y,r)&=\bar{e}\bigl[\eta_{\lambda}(r)\beta(x,y)\bigr]\overline{f(-x,-y,r)}.
\end{align}

 Meanwhile, the multiplicative unitary operator defines the {\em comultiplication\/}
$\Delta:A\to M(A\otimes A)$.  For $a\in A$, we have:
$$
\Delta a=U_A(a\otimes1){U_A}^*.
$$
At the level of functions in ${\mathcal A}$, the equation $\Delta(L_f)
=(L\otimes L)_{\Delta f}$ gives us the following:
\begin{align}\label{(Acomult)}
&\Delta f(x,y,r;x',y',r')  \notag \\
&=\int f(x',y',r+r')\bar{e}\bigl[\tilde{p}\cdot(e^{\lambda r'}x'-x)+\tilde{q}
\cdot(e^{\lambda r'}y'-y)\bigr]\,d\tilde{p}d\tilde{q},
\end{align}
which is a Schwartz function having compact support in $r$ and $r'$.

 There is also the map $S:A\to A$, defined by $S(a)=\hat{J}a^*\hat{J}$, where
$\hat{J}$ is the following involutive operator on ${\mathcal H}$:
$$
\hat{J}\xi(x,y,r)=(e^{\lambda r})^n\overline{\xi(e^{\lambda r}x,e^{\lambda r}y,-r)}.
$$
See Proposition 2.4 of \cite{BJKppha}.  Then $S(L_f)=L_{S(f)}$ gives us the
following:
\begin{equation}\label{(Aantipode)}
\bigl(S(f)\bigr)(x,y,r)=(e^{2\lambda r})^n\bar{e}\bigl[\eta_{\lambda}(r)\beta(x,y)
\bigr]f(-e^{\lambda r}x,-e^{\lambda r}y,-r),
\end{equation}
at the level of functions in ${\mathcal A}$.

 Turning our focus to the dual object of $(A,\Delta)$, we now consider the
``right slices'' of $U_A$.  That is, the $C^*$-algebra $\hat{A}$ is generated
by the operators $\rho(\omega)=(\operatorname{id}\otimes\omega)(U_A)$, for
$\omega\in{\mathcal B}({\mathcal H})_*$:
$$
\hat{A}=\overline{\bigl\{(\operatorname{id}\otimes\omega)(U_A):\omega\in
{\mathcal B}({\mathcal H})_*\bigr\}}^{\|\ \|}.
$$
The comultiplication $\hat{\Delta}:\hat{A}\to M(\hat{A}\otimes\hat{A})$ is
given by $\hat{\Delta}b={U_A}^*(1\otimes b)U_A$.

 As above, there is an alternative characterization of $\hat{A}$.  For this,
consider $\hat{\mathcal A}$, which is also the space of Schwartz functions
in the $(x,y,r)$ variables having compact support in $r$.  Define the ``regular
representation'' $\rho$ of $\hat{\mathcal A}$ on ${\mathcal B}({\mathcal H})$,
given by
$$
(\rho_{\phi}\xi)(x,y,r):=\int (e^{\lambda\tilde{r}})^n\phi(x,y,\tilde{r})
\xi(e^{\lambda\tilde{r}}x,e^{\lambda\tilde{r}}y,r-\tilde{r})\,d\tilde{r}.
$$
We have seen in \cite{BJKqhg} that $\hat{A}=\overline{\rho(\hat
{\mathcal A})}^{\|\ \|}$.  As before, we can regard the functions
$\phi\in\hat{\mathcal A}$ same as the operators $\rho_{\phi}$, and
$\hat{\mathcal A}$ is considered as a dense ${}^*$-subalgebra of
$\hat{A}$.  On $\hat{\mathcal A}$, the multiplication and the involution
take the following form (via $\rho_{\phi\times\psi}=\rho_{\phi}\rho_{\psi}$
and $\rho_{\phi^*}=(\rho_{\phi})^*$):
\begin{align}\label{(hatAops)}
(\phi\times_{\hat{A}}\psi)(x,y,r)&=\int\phi(x,y,\tilde{r})\psi(e^{\lambda
\tilde{r}}x,e^{\lambda\tilde{r}}y,r-\tilde{r})\,d\tilde{r}.  \notag \\
\phi^*(x,y,r)&=\overline{\phi(e^{\lambda r}x,e^{\lambda r}y,-r)}.
\end{align}

 Meanwhile, we have the following description of the comultiplication,
obtained by the equation $\hat{\Delta}(\rho_{\phi})=(\rho\otimes\rho)_{\hat
{\Delta}\phi}$:
\begin{align}\label{(hatAcomult)}
&\hat{\Delta}\phi(x,y,r;x',y',r')  \notag \\
&=\int\phi(x+x',y+y',\tilde{r})e\bigl[\eta_{\lambda}(\tilde{r})\beta(x,y')
\bigr]e\bigl[\tilde{r}(z+z')\bigr]\bar{e}[zr+z'r']\,d\tilde{r}dzdz'.
\end{align}
The antipode $\hat{S}:\hat{A}\to\hat{A}$ is given by $\hat{S}=Jb^*J$, where
$J$ is the operator on ${\mathcal H}$ defined by
$$
J\xi(x,y,r)=\bar{e}\bigl[\eta_{\lambda}(r)\beta(x,y)\bigr]\overline
{\xi(-x,-y,r)}.
$$
Then at the level of functions in $\hat{\mathcal A}$, the expression $\hat{S}
(\rho_{\phi})=\rho_{\hat{S}(\phi)}$ gives us the following:
\begin{equation}\label{(hatAantipode)}
\bigl(\hat{S}(\phi)\bigr)(x,y,r)=\bar{e}\bigl[\eta_{\lambda}(r)\beta(x,y)
\bigr]\phi(-e^{\lambda r}x,-e^{\lambda r}y,-r).
\end{equation}

 We have further shown in our previous papers that the two $C^*$-bialgebras
$(A,\Delta)$ and $(\hat{A},\hat{\Delta})$ are indeed examples of non-compact,
$C^*$-algebraic quantum groups (together with the necessary ingredients like
Haar weights).  See \cite{BJKppha} and \cite{BJKqhg}.  They are mutually dual
objects in the framework of locally compact quantum groups.  Since the square
of the antipode map is identity for both of them (which can be easily seen
from the definitions of $S$ and $\hat{S}$ given above), they are cases of
{\em Kac $C^*$-algebras\/} (as in \cite{Va}).

 Unlike in the purely algebraic or finite-dimensional setting, no natural
dual pairing exists between $A$ and $\hat{A}$. However, at least at the level
of the dense subalgebras ${\mathcal A}$ and $\hat {\mathcal A}$, there does
exist a suitable dual pairing, defined as follows:
\begin{equation}\label{(dualpairing)}
\langle f,\phi\rangle=\int f(x,y,r)\phi(e^{\lambda r}x,e^{\lambda r}y,-r)\,dxdydr,
\end{equation}
for $f(=L_f)\in{\mathcal A}$ and $\phi(=\rho_{\phi})\in\hat{\mathcal A}$.  See
Proposition 3.1 of \cite{BJKqhg}, which is just an immediate consequence of
Definition 1.3 of \cite{BS}.  As we have shown in Proposition 3.1 of \cite
{BJKqhg}, this dual pairing satisfies all the necessary properties for it to
be considered as the correct analog of the pairing in the (purely algebraic)
framework of Hopf algebras.

 Since we are planning to construct the quantum double, we also need to
clarify the ``opposite'' and ``co-opposite'' versions of $(A,\Delta)$
and $(\hat{A},\hat{\Delta})$, again in the $C^*$-algebra framework.
For this, it is quite helpful to know that we can form a {\em Kac system\/}
(in the sense of Baaj and Skandalis) from our multiplicative unitary operator
$U_A$.  The following observation was made in section 3 of \cite{BJKqhg}:

\begin{prop}\label{jKacsystem}
Let $j\in{\mathcal B}({\mathcal H})$ be defined by $j=\hat{J}J=J\hat{J}$,
where $J$ and $\hat{J}$ are the anti-unitary operators appeared in the
definitions of the antipode maps.  Then $j$ is an (involutive) unitary operator
given by
$$
j\xi(x,y,r)=(e^{\lambda r})^n\bar{e}\bigl[\eta_{\lambda}(r)\beta(x,y)\bigr]
\xi(-e^{\lambda r}x,-e^{\lambda r}y,-r).
$$
Moreover, the triple $({\mathcal H},U_A,j)$ forms a ``Kac system'' (as in
section 6 of \cite{BS}).  In particular, the following unitary operators
are all multiplicative:
\begin{align}
&U_A\in M(\hat{A}\otimes A), \notag \\
&\widehat{U_A}=\Sigma(j\otimes1)U_A(j\otimes1)\Sigma\in M(A\otimes
\hat{A}^{\operatorname{op}}),  \notag \\
&\widetilde{U_A}=(j\otimes j)\widehat{U_A}(j\otimes j)=(j\otimes1)(\Sigma
U_A\Sigma)(j\otimes1)\in M(A^{\operatorname{op}}\otimes\hat{A}), \notag \\
&\widehat{\widehat{U_A}}=\widetilde{\widetilde{U_A}}=(j\otimes j)U_A
(j\otimes j)\in M(\hat{A}^{\operatorname{op}}\otimes A^{\operatorname{op}}),
\notag
\end{align}
where $\Sigma$ denotes the flip.
\end{prop}

\begin{rem}
This proposition is actually a consequence of the fact that $(A,\Delta)$
and $(\hat{A},\hat{\Delta})$ are Kac $C^*$-algebras.  The results may be
checked by a direct computation.  See Proposition 3.2 of \cite{BJKqhg},
and also 6.11(d) of \cite{BS}.
\end{rem}

 Using the multiplicative unitary operators $U_A$ and its variations
obtained in the above, we can define several different versions of
the quantum Heisenberg group algebra and the quantum Heisenberg group,
in the form of $(A^{\operatorname{op}},{\Delta})$, $(A,{\Delta}^{\operatorname
{cop}})$, $(A^{\operatorname{op}},{\Delta}^{\operatorname{cop}})$, as well as
$(\hat{A}^{\operatorname{op}},\hat{\Delta})$, $(\hat{A},\hat{\Delta}^{\operatorname
{cop}})$, $(\hat{A}^{\operatorname{op}},\hat{\Delta}^{\operatorname{cop}})$.

 For instance, $(\hat{A}^{\operatorname{op}},\hat{\Delta})$ is determined by
the multiplicative unitary operator
$X=\Sigma\widehat{U_A}^*\Sigma$.  To be more precise, we have:
$$
\hat{A}^{\operatorname{op}}=\overline{\bigl\{(\operatorname{id}\otimes\omega)
(X):\omega\in{\mathcal B}({\mathcal H})_*\bigr\}}^{\|\ \|}=\overline{\lambda
(\hat{\mathcal A})}^{\|\ \|}\bigl(\subseteq{\mathcal B}({\mathcal H})\bigr),
$$
where $\lambda:\hat{\mathcal A}\to{\mathcal B}({\mathcal H})$ is defined by
$$
(\lambda_{\phi}\xi)(x,y,r):=\int\phi(e^{\lambda\tilde{r}}x,e^{\lambda\tilde{r}}y,
r-\tilde{r})\xi(x,y,\tilde{r})\,d\tilde{r}.
$$
Notice that $\lambda_{\phi}\lambda_{\psi}=\lambda_{\psi\times\phi}$,
implementing the opposite multiplication.  The comultiplication,
given by $\hat{A}^{\operatorname{op}}\ni b\mapsto X^*(1\otimes b)X
\in M(\hat{A}^{\operatorname{op}}\otimes\hat{A}^{\operatorname{op}})$,
stays the same at the function level: That is, $X^*(1\otimes\lambda_{\phi})X
=(\lambda\otimes\lambda)_{\hat{\Delta}\phi}$, where $\hat{\Delta}\phi$ is
same as in \eqref{(hatAcomult)}.  The antipode also stays the same: $\hat{S}
(\lambda_{\phi})=\lambda_{\hat{S}(\phi)}$,  as in equation \eqref{(hatAantipode)}.

 We only gave here one possible description of $(\hat{A}^{\operatorname{op}},
\hat{\Delta})$, since $\hat{A}^{\operatorname{op}}$ is the one we immediately
need for the definition of our quantum double.  But See Proposition 3.5 of
\cite{BJKqhg} for the others.

\section{The quantum double}

 Since we know the expressions for various operations on ${\mathcal A}$ and
$\hat{\mathcal A}$ (the equations \eqref{(Aops)}, \eqref{(Acomult)}, \eqref
{(Aantipode)}, and \eqref{(hatAops)}, \eqref{(hatAcomult)}, \eqref{(hatAantipode)}),
as well as the expression for the dual pairing between them given by equation
\eqref{(dualpairing)}, we can use Definition \ref{qdouble} to write down the
product on the quantum double, at the function level:
\begin{align}\label{(doubleprod)}
&\bigl((\phi\otimes f)\times(\psi\otimes g)\bigr)(x,y,r;x',y',r') \notag \\
&=\int\phi(e^{\lambda\tilde{r}}x,e^{\lambda\tilde{r}}y,r-\tilde{r})\psi(x-e^{\lambda
(r'-\tilde{r})}\tilde{x}+e^{-\lambda\tilde{r}}\tilde{x},y-e^{\lambda(r'-\tilde{r})}
\tilde{y}+e^{-\lambda\tilde{r}}\tilde{y},\tilde{r})  \notag \\
&\qquad \bar{e}\bigl[\eta_{\lambda}(r')\beta(e^{-\lambda\tilde{r}}\tilde{x},y')\bigr]
e\bigl[\eta_{\lambda}(\tilde{r})\beta(x,e^{-\lambda\tilde{r}}\tilde{y})\bigr]
\bar{e}\bigl[\eta_{\lambda}(\tilde{r})\beta(e^{\lambda(r'-\tilde{r})}\tilde{x},y)\bigr]
\notag \\
&\qquad e[\eta_{\lambda}(r')\beta(\tilde{x},\tilde{y})\bigr]e\bigl[\eta_{\lambda}
(\tilde{r})\beta(e^{-\lambda\tilde{r}}\tilde{x},e^{-\lambda\tilde{r}}\tilde{y})\bigr]
\bar{e}\bigl[\eta_{\lambda}(\tilde{r})\beta(e^{\lambda(r'-\tilde{r})}\tilde{x},
e^{-\lambda\tilde{r}}\tilde{y})\bigr]  \notag \\
&\qquad f(\tilde{x},\tilde{y},r')g(x'-e^{-\lambda\tilde{r}}\tilde{x},
y'-e^{-\lambda\tilde{r}}\tilde{y},r')\,d\tilde{x}d\tilde{y}d\tilde{r}.
\end{align}
Here, $\phi,\psi\in\hat{\mathcal A}$ and $f,g\in{\mathcal A}$.  Computation is
rather long, but not really difficult.

 However, for us to be able to define $D(A)$ properly at the $C^*$-algebra
level, it is again best to work with the multiplicative unitary operators.
Since we wish to construct an object that will be considered as containing
$(A,\Delta)$ and its ``opposite dual'' $({\hat{A}}^{\operatorname{op}},
\hat{\Delta})$, with some actions involved, let us define the following
unitary operator:
\begin{equation}\label{(V_D)}
V_D=Z_{12}Y_{24}Z_{12}^*X_{13}\in{\mathcal B}({\mathcal H}\otimes{\mathcal H}
\otimes{\mathcal H}\otimes{\mathcal H}).
\end{equation}
Here $X=\Sigma\widehat{U_A}^*\Sigma$ is as defined in the previous section,
$Y=\Sigma X^*\Sigma=\widehat {U_A}$, while $Z=Y^*\widehat{\widehat{Y}}=\widehat
{\widehat{Y}}Y^*$ (It is known that $Y\in M(A\otimes\hat{A}^{\operatorname{op}})$
and $\widehat{\widehat{Y}}\in M(A^{\operatorname{op}}\otimes\hat{A})$, which
can be seen in Corollary of Proposition 3.5 of \cite{BJKqhg}.  See also
Proposition~\ref{jKacsystem}.).  The leg numbering notation is as before.

 Main motivation for our choice comes from section 8 of \cite{BS}, and
the formulation is essentially equivalent to the ones given in \cite{Ya},
\cite{BjVa} (though slightly different). Roughly speaking, the operator
$X$ gives $({\hat{A}}^{\operatorname{op}},\hat{\Delta})$ (as we saw in
Section 2), the operator $Y$ gives $(A,\Delta)$ (as in Proposition 3.5 of
\cite{BJKqhg}), and the operator $Z$ enables us to encode the generalized
coadjoint actions.  See Proposition \ref{Zoperator} below, which comes after
the following short lemma:

\begin{lem}\label{ZY^*}
For any $a\in A$ and $b\in\hat{A}^{\operatorname{op}}$, we have:
$$
Z(a\otimes b)Z^*=Y^*(a\otimes b)Y.
$$
\end{lem}

\begin{proof}
The proof easily follows from the fact that $Z=Y^*\widehat{\widehat{Y}}
=\widehat{\widehat{Y}}Y^*$, while $Y\in M(A\otimes\hat{A}^{\operatorname{op}})$
and $\widehat{\widehat{Y}}\in M(A^{\operatorname{op}}\otimes\hat{A})$.  Actually,
the result will still hold if $a\in M(A)$ and $b\in M(\hat{A}^{\operatorname{op}})$.
\end{proof}

\begin{prop}\label{Zoperator}
Let $Z=Y^*\widehat{\widehat{Y}}=\widehat{\widehat{Y}}Y^*\in{\mathcal B}
({\mathcal H}\otimes{\mathcal H})$ be the operator defined above.  Explicitly,
\begin{align}
Z\xi(x,y,r;x',y',r')
&=\bar{e}\bigl[\eta_{\lambda}(r)\beta(e^{\lambda r'}x',y-e^{\lambda r'}y')\bigr]
e\bigl[\eta_{\lambda}(r)\beta(x-e^{\lambda r'}x',y')\bigr] 
\notag \\
&\quad(e^{\lambda r})^n\,\xi(x-e^{\lambda r'}x'+x',y-e^{\lambda r'}y'+y',r;
e^{\lambda r}x',e^{\lambda r}y',r').  \notag
\end{align}
Let $\alpha:A\to M(\hat{A}^{\operatorname{op}}\otimes A)$ and $\alpha':
\hat{A}^{\operatorname{op}}\to M(\hat{A}^{\operatorname{op}}\otimes A)$
be defined by
$$
\alpha(a):=\Sigma Z^*(a\otimes1)Z\Sigma,\qquad {\text {and }}\qquad \alpha'(b):=
\Sigma Z^*(1\otimes b)Z\Sigma.
$$
Then $\alpha$ is a left coaction of $(A,\Delta)$ on the $C^*$-algebra
$\hat{A}^{\operatorname{op}}$, while $\alpha'$ is a right coaction of
$(\hat{A}^{\operatorname{op}},\hat{\Delta})$ on $A$.  That is, the maps
$\alpha$ and $\alpha'$ are non-degenerate ${}^*$-homomorphisms such that:
$$
(\hat{\Delta}\otimes\operatorname{id})\alpha=(\operatorname{id}\otimes\alpha)
\alpha\qquad {\text {and}}\qquad(\operatorname{id}\otimes\Delta)\alpha'
=(\alpha'\otimes\operatorname{id})\alpha'.
$$
\end{prop}

\begin{proof}
Let $a\in A$.  Then we have:
\begin{align}
(\operatorname{id}\otimes\alpha)\alpha(a)&=\Sigma_{23}Z_{23}^*\Sigma_{12}
Z_{12}^*(a\otimes1\otimes1)Z_{12}\Sigma_{12}Z_{23}\Sigma_{23}  \notag \\
&=\Sigma_{23}Y_{23}\Sigma_{12}Y_{12}(a\otimes1\otimes1)Y_{12}^*
\Sigma_{12}Y_{23}^*\Sigma_{23}  \notag \\
&=X_{23}^*\Sigma_{23}X_{12}^*\Sigma_{12}(a\otimes1\otimes1)\Sigma_{12}
X_{12}\Sigma_{23}X_{23}  \notag \\
&=X_{23}^*X_{13}^*(1\otimes1\otimes a)X_{13}X_{23}.  \notag
\end{align}
In the second equality, we are using Lemma \ref{ZY^*}.  And in the third
equality, we are using that $Y=\Sigma X^*\Sigma$.  On the other hand,
remembering that $\hat{\Delta}(b)=X^*(1\otimes b)X$, for $b\in
\hat{A}^{\operatorname{op}}$, we have:
\begin{align}
(\hat{\Delta}\otimes\operatorname{id})\alpha(a)&=X_{12}^*\Sigma_{23}Z_{23}^*
(1\otimes a\otimes1)Z_{23}\Sigma_{23}X_{12}  \notag \\
&=X_{12}^*X_{23}^*\Sigma_{23}(1\otimes a\otimes1)\Sigma_{23}X_{23}X_{12}
\notag \\
&=X_{12}^*X_{23}^*(1\otimes1\otimes a)X_{23}X_{12}  \notag \\
&=X_{23}^*X_{13}^*X_{12}^*(1\otimes1\otimes a)X_{12}X_{13}X_{23}  \notag \\
&=X_{23}^*X_{13}^*(1\otimes1\otimes a)X_{13}X_{23}.  \notag
\end{align}
We again used Lemma \ref{ZY^*} and $Y=\Sigma X^*\Sigma$ in the second equality
above.  In the fourth equality, the multiplicativity of $X$ (satisfying the
pentagon equation: $X_{12}X_{13}X_{23}=X_{23}X_{12}$) was used.  In this way,
we show that: $(\hat{\Delta}\otimes\operatorname{id})\alpha
=(\operatorname{id}\otimes\alpha)\alpha$.

The condition for $\alpha'$ is similarly proved, using that $\Delta a
=Y^*(1\otimes a)Y$, for $a\in A$ (This is noted in Proposition 3.5 of
\cite{BJKqhg}.).
\end{proof}

\begin{rem}
The coactions $\alpha$ and $\alpha'$ are essentially ``coadjoint
coactions'', which (dually) correspond to the ``coadjoint actions''
given in Definition~\ref{qdouble}.  Indeed, at least at the level
of functions in ${\mathcal A}$ and $\hat{\mathcal A}$, it is possible
to show that:
$$
\bigl\langle\alpha(a),f\otimes\phi\bigr\rangle=\langle a,f\triangleright
\phi\rangle \qquad {\text {and}}\qquad \bigl\langle\alpha'(b),f\otimes\phi
\bigr\rangle=\langle b,f\triangleleft\phi\rangle,
$$
where $a,f\in{\mathcal A}$ and $b,\phi\in\hat{\mathcal A}$, while
$\langle\ ,\ \rangle$ is the dual pairing as given in equation
\eqref{(dualpairing)}.
\end{rem}

 In the ensuing paragraphs, we will show that the operator $V_D$ as defined
in equation \eqref{(V_D)} is actually multiplicative, and make our case that
the $C^*$-bialgebra it determines is exactly the quantum double $D(A)$ we
are looking for.  In particular, we will see that the $C^*$-algebra contains
as a dense subalgebra $\hat{\mathcal A}\otimes{\mathcal A}$, with its product
defined in equation \eqref{(doubleprod)}. 

 The multiplicativity of $V_D$ (i.\,e. satisfying the ``pentagon equation'')
could be shown directly, but the computation will be rather long and tedious
due to the fact that we have to work with 18 variables.  So we present here
an alternative way, where the crucial point is that the operators $X$ and $Y$
give rise to a ``matched pair'' (See Definition 8.13 of \cite{BS}) of Kac systems.

\begin{lem}\label{lemmaV}
Let the notation be as above.
\begin{enumerate}
\item The triples $({\mathcal H},X,j)$ and $({\mathcal H},Y,j)$, together with
the operator $Z$, form a matched pair of Kac systems.
\item The operator $V:=Z_{12}^*X_{13}Z_{12}Y_{24}$ is multiplicative.
\item $Z_{34}V=Z_{34}Z_{12}^*X_{13}Z_{12}Y_{24}=Y_{24}Z_{12}^*X_{13}Z_{12}
Z_{34}$.
\end{enumerate}
\end{lem}

\begin{proof}
(1) is the result of Theorem 8.17 of \cite{BS}, from which the multiplicativity
of $V$ follows (By Definition 8.15 of \cite{BS}, the operator $V$ determines
the ``$Z$-tensor product'' of the matched pair.).  See also our Proposition
\ref{tauinversion} and its Corollary in Section 4 below.  Meanwhile, by
Proposition 8.10 of \cite{BS}, condition (2) is equivalent to condition (3)
(We can also check (3) directly from the definitions.).
\end{proof}

\begin{cor}
The unitary operator $V_D$ defined in \eqref{(V_D)} is multiplicative.
\end{cor}

\begin{proof}
By (3) of the previous lemma, we have:
$$
Z_{12}Z_{34}VZ_{34}^*Z_{12}^*=Z_{12}Y_{24}Z_{12}^*X_{13}=V_D.
$$
Re-writing this expression as $V_D=(Z\otimes Z)V(Z^*\otimes Z^*)$
and noting that $Z$ is unitary, we see easily that $V_D$ is also
a multiplicative unitary operator (by being unitarily equivalent
to $V$).
\end{proof}

 By the general theory of multiplicative unitary operators \cite{BS}, \cite{Wr7},
the operator $V_D$ will let us define a $C^*$-bialgebra, on which we build
the necessary ingredients for it to become a locally compact quantum group.
Specifically, let us consider the $C^*$-bialgebra $(A_D,\Delta_D)$, which is
generated by the ``right slices'' of $V_D$, as follows:

\begin{defn}\label{defnA_D}
Let $A_D$ be the $C^*$-algebra contained in ${\mathcal B}({\mathcal H}
\otimes{\mathcal H})$, defined by
$$
A_D=\overline{\bigl\{(\operatorname{id}\otimes\operatorname{id}\otimes
\Omega)(V_D):\Omega\in{\mathcal B}({\mathcal H}\otimes{\mathcal H})_*
\bigr\}}^{\|\ \|}.
$$
For a typical element $x\in A_D$, define $\Delta_D(x)$ by
$$
\Delta_D(x):={V_D}^*(1\otimes1\otimes x)V_D.
$$
In this way, we obtain the {\em comultiplication\/} $\Delta_D:A_D\to
M(A_D\otimes A_D)$, which is a non-degenerate $C^*$-homomorphism satisfying
the coassociativity condition: $(\Delta_D\otimes\operatorname{id})\Delta_D
=(\operatorname{id}\otimes\Delta_D)\Delta_D$.  Moreover, $\Delta_D(A_D)
(A_D\otimes1)$ and $\Delta_D(A_D)(1\otimes A_D)$ are dense subsets in
$A_D\otimes A_D$.
\end{defn}

 For the last statement (the density conditions), see Theorem 1.5 and
section 5 of \cite{Wr7} (It is a non-trivial result.).  Our goal now is
to show that $(A_D,\Delta_D)$ is exactly the quantum double $D(A)$, analogous
to Definition \ref{qdouble}.  Let us first give a more concrete $C^*$-algebraic
realization of $A_D$.

\begin{prop}\label{A_D}
Let $\pi:A\to{\mathcal B}({\mathcal H}\otimes{\mathcal H})$ and $\pi':
\hat{A}^{\operatorname{op}}\to{\mathcal B}({\mathcal H}\otimes{\mathcal H})$
be defined by
$$
\pi(a):=Z(1\otimes a)Z^*\qquad {\text {and}}\qquad \pi'(b):=b\otimes1.
$$
Then $A_D$ is the $C^*$-algebra generated by the operators $\pi(a)\pi'(b)$,
for $a\in A$, $b\in\hat{A}^{\operatorname{op}}$.
\end{prop}

\begin{proof}
For $\omega,\omega'\in{\mathcal B}({\mathcal H})_*$, we have:
\begin{align}
(\operatorname{id}\otimes\operatorname{id}\otimes\omega\otimes\omega')(V_D)
&=(\operatorname{id}\otimes\operatorname{id}\otimes\omega\otimes\omega')
(Z_{12}Y_{24}Z_{12}^*X_{13}) \notag \\
&=Z\bigl[1\otimes(\operatorname{id}\otimes\omega')(Y)\bigr]Z^*
\bigl[(\operatorname{id}\otimes\omega)(X)\otimes1\bigr]  \notag \\
&=\pi(a)\pi'(b),  \notag
\end{align}
where $a=(\operatorname{id}\otimes\omega')(Y)$ and $b=(\operatorname{id}\otimes
\omega)(X)$.  This is valid, because by Proposition 3.5\,(2) of \cite{BJKqhg}
and by section 6 (Appendix) of \cite{BJKppha}, we have:
$$
a=(\operatorname{id}\otimes\omega')(Y)=(\operatorname{id}\otimes\omega')(\widehat{U_A})
\in A
$$
and
$$
b=(\operatorname{id}\otimes\omega)(X)=(\operatorname{id}\otimes\omega)(\Sigma
\widehat{U_A}^*\Sigma)\in\hat{A}^{\operatorname{op}}.
$$
In fact, the operators $(\operatorname{id}\otimes\omega')(Y)$, $\omega'\in{\mathcal B}
({\mathcal H})_*$, generate $A$, while the operators $(\operatorname{id}\otimes\omega)
(X)$, $\omega\in{\mathcal B}({\mathcal H})_*$, generate $\hat{A}^{\operatorname{op}}$.

Since the operators $(\operatorname{id}\otimes\operatorname{id}\otimes\omega\otimes
\omega')(V_D)$ generate $A_D$ (Definition~\ref{defnA_D}), the result of the proposition
follows.
\end{proof}

\begin{cor}
The maps $\pi$ and $\pi'$ are $C^*$-algebra homomorphisms.  Namely,
$$
\pi:A\to M(A_D)\qquad {\text {and}}\qquad \pi':\hat{A}^{\operatorname{op}}\to
M(A_D).
$$
\end{cor}

\begin{rem}
The corollary is obvious from the definitions of $\pi$ and $\pi'$.  Here, $M(A_D)$
denotes the multiplier algebra of $A_D$.  Later, when we clarify the co-algebra
structure on $A_D$, they will actually become $C^*$-bialgebra homomorphisms.
\end{rem}

\begin{prop}\label{Pi}
Let $\Pi:\hat{\mathcal A}\otimes{\mathcal A}\to{\mathcal B}({\mathcal H}\otimes
{\mathcal H})$ be defined by
$$
\Pi(\phi\otimes f):=\pi'(\lambda_{\phi})\pi(L_f),\qquad {\text {for $\phi\in
\hat{\mathcal A}, f\in{\mathcal A}$.}}
$$
Then $A_D=\overline{\Pi(\hat{\mathcal A}\otimes{\mathcal A})}^{\|\ \|}$.
\end{prop}

\begin{proof}
We know from Section 2 that $A=\overline{L({\mathcal A})}^{\|\ \|}$ and
$\hat{A}^{\operatorname{op}}=\overline{\lambda(\hat{\mathcal A})}^{\|\ \|}$.
So this result is an immediate consequence of Proposition \ref{A_D}, with the
aid of the fact that the $C^*$-algebras are closed under involution.
\end{proof}

 We observe that $\Pi$ determines a multiplication on $\hat{\mathcal A}\otimes
{\mathcal A}$, given by
$$
\Pi(\phi\otimes f)\Pi(\psi\otimes g)=\Pi\bigl((\phi\otimes f)\times(\psi\otimes g)
\bigr),
$$
making it a (dense) subalgebra of $A_D$.  It turns out that the product obtained
in this way exactly coincides with the one given in \eqref{(doubleprod)}.

\begin{prop}\label{prodA_D}
Let $\hat{\mathcal A}\otimes{\mathcal A}$ be given the multiplication, as in
\eqref{(doubleprod)}.  Then we have:
$$
\Pi(\phi\otimes f)\Pi(\psi\otimes g)=\Pi\bigl((\phi\otimes f)\times(\psi\otimes g)\bigr),
$$
for $\phi,\psi\in\hat{\mathcal A}$ and $f,g\in{\mathcal A}$.
\end{prop}

\begin{proof}
For $\phi\in\hat{\mathcal A}$ and $\xi\in{\mathcal H}\otimes{\mathcal H}$,
\begin{align}
\pi'(\lambda_{\phi})\xi(x,y,r;x',y',r')&=(\lambda_{\phi}\otimes1)\xi(x,y,r;x',y',r')
\notag \\
&=\int\phi(e^{\lambda\tilde{r}}x,e^{\lambda\tilde{r}}y,r-\tilde{r})\xi(x,y,\tilde{r};
x',y',r')\,d\tilde{r}.  \notag
\end{align}
Whereas for $f\in{\mathcal A}$ and $\xi\in{\mathcal H}\otimes{\mathcal H}$,
\begin{align}
&\pi(L_f)\xi(x,y,r;x',y',r')=Z(1\otimes L_f)Z^*\xi(x,y,r;x',y',r')  \notag \\
&=(e^{\lambda r})^n\,\bar{e}\bigl[\eta_{\lambda}(r)\beta(e^{\lambda r'}x',
y-e^{\lambda r'}y')\bigr]e\bigl[\eta_{\lambda}(r)\beta(x-e^{\lambda r'}x',y')\bigr] 
\notag \\
&\quad(1\otimes L_f)Z^*\xi(x-e^{\lambda r'}x'+x',y-e^{\lambda r'}y'+y',r;
e^{\lambda r}x',e^{\lambda r}y',r')  \notag \\
&=\int(e^{\lambda r})^n\,\bar{e}\bigl[\eta_{\lambda}(r)\beta(e^{\lambda r'}x',
y-e^{\lambda r'}y')\bigr]e\bigl[\eta_{\lambda}(r)\beta(x-e^{\lambda r'}x',y')\bigr] 
\notag \\
&\qquad f(\tilde{x},\tilde{y},r')\bar{e}\bigl[\eta_{\lambda}(r')\beta(\tilde{x},
e^{\lambda r}y'-\tilde{y})\bigr] \notag \\
&\qquad Z^*\xi(x-e^{\lambda r'}x'+x',y-e^{\lambda r'}y'+y',r;
e^{\lambda r}x'-\tilde{x},e^{\lambda r}y'-\tilde{y},r')\,d\tilde{x}d\tilde{y}
\notag \\
&=(\cdots)  \notag \\
&=\int f(\tilde{x},\tilde{y},r')\bar{e}\bigl[\eta_{\lambda}(r')\beta(e^{-\lambda r}
\tilde{x},y')\bigr]e\bigl[\eta_{\lambda}(r)\beta(x,e^{-\lambda r}\tilde{y})\bigr]
\bar{e}\bigl[\eta_{\lambda}(r)\beta(e^{\lambda r'-\lambda r}\tilde{x},y)\bigr]
\notag\\
&\qquad e\bigl[\eta_{\lambda}(r')\beta(\tilde{x},\tilde{y})\bigr]
e\bigl[\eta_{\lambda}(r)\beta(e^{-\lambda r}\tilde{x},e^{-\lambda r}\tilde{y})\bigr]
\bar{e}\bigl[\eta_{\lambda}(r)\beta(e^{\lambda r'-\lambda r}\tilde{x},e^{-\lambda r}
\tilde{y})\bigr]  \notag \\
&\qquad\xi(x-e^{\lambda r'-\lambda r}\tilde{x}+e^{-\lambda r}\tilde{x},
y-e^{\lambda r'-\lambda r}\tilde{y}+e^{-\lambda r}\tilde{y},r;x'-e^{-\lambda r}\tilde{x},
y'-e^{-\lambda r}\tilde{y},r')\,d\tilde{x}d\tilde{y}.  \notag
\end{align}
So we have:
\begin{align}
&\Pi(\phi\otimes f)\xi(x,y,r;x',y',r')=\pi'(\lambda_{\phi})\pi(L_f)\xi(x,y,r;x',y',r')
\notag \\
&=\int\phi(e^{\lambda\tilde{r}}x,e^{\lambda\tilde{r}}y,r-\tilde{r})f(\tilde{x},\tilde{y},r')
\bar{e}\bigl[\eta_{\lambda}(r')\beta(e^{-\lambda\tilde{r}}\tilde{x},y')\bigr]  \notag \\
&\qquad e\bigl[\eta_{\lambda}(\tilde{r})\beta(x,e^{-\lambda\tilde{r}}\tilde{y})\bigr]
\bar{e}\bigl[\eta_{\lambda}(\tilde{r})\beta(e^{\lambda r'-\lambda\tilde{r}}\tilde{x},y)\bigr]
\notag \\
&\qquad e[\eta_{\lambda}(r')\beta(\tilde{x},\tilde{y})\bigr]e\bigl[\eta_{\lambda}
(\tilde{r})\beta(e^{-\lambda\tilde{r}}\tilde{x},e^{-\lambda\tilde{r}}\tilde{y})\bigr]
\bar{e}\bigl[\eta_{\lambda}(\tilde{r})\beta(e^{\lambda r'-\lambda\tilde{r}}\tilde{x},
e^{-\lambda\tilde{r}}\tilde{y})\bigr]  \notag \\
&\qquad \xi(x-e^{\lambda r'-\lambda\tilde{r}}\tilde{x}+e^{-\lambda\tilde{r}}\tilde{x},
y-e^{\lambda r'-\lambda\tilde{r}}\tilde{y}+e^{-\lambda\tilde{r}}\tilde{y},\tilde{r};
\notag \\
&\qquad\quad x'-e^{-\lambda\tilde{r}}\tilde{x},y'-e^{-\lambda\tilde{r}}\tilde{y},
r')\,d\tilde{x}d\tilde{y}d\tilde{r}.  \notag
\end{align}
Using this and noting its resemblance to equation \eqref{(doubleprod)},
it is not difficult to show that:
$$
\Pi(\phi\otimes f)\Pi(\psi\otimes g)\xi=\Pi\bigl((\phi\otimes f)\times(\psi\otimes g)
\bigr)\xi,\qquad {\text {for any $\xi\in{\mathcal H}\otimes{\mathcal H}$.}}
$$
\end{proof}

 The involution on $A_D$ is inherited from that of ${\mathcal B}({\mathcal H}\otimes
{\mathcal H})$.  At the level of the subalgebra $\hat{\mathcal A}\otimes{\mathcal A}$,
it takes the following form:
\begin{align}\label{(doubleinvol)}
&(\phi\otimes f)^*(x,y,r;x',y',r') \notag \\
&=(e^{2\lambda r})^n \overline{f(-e^{\lambda r}x',-e^{\lambda r}y',r')}
e\bigl[\eta_{\lambda}(r)\beta(x,y')\bigr]\bar{e}\bigl
[\eta_{\lambda}(r)\beta(e^{\lambda r'}x',y)\bigr]  \notag \\
&\quad\bar{e}\bigl[\eta_{\lambda}(r'-r)\beta(e^{\lambda r}x',e^{\lambda r}y')\bigr]
\bar{e}\bigl[\eta_{\lambda}(r)\beta (e^{\lambda r'}x',y')\bigr]  \notag \\
&\quad\overline{\phi(e^{\lambda r}x-e^{\lambda r'+\lambda r}x'+e^{\lambda r}x',
e^{\lambda r}y-e^{\lambda r'+\lambda r}y'+e^{\lambda r}y',-r)}.
\end{align}
To be more precise, the definition of $(\phi\otimes f)^*\in\hat{\mathcal A}\otimes
{\mathcal A}$ above has been chosen so that we have: $\Pi\bigl((\phi\otimes f)^*\bigr)
=\bigl[\Pi(\phi\otimes f)\bigr]^*=\pi(L_f)^*\pi'(\lambda_{\phi})^*$.

 Next, let us turn our attention to the co-algebra structure on $A_D$.  As in the
previous proposition, we will see that at the level of functions in $\hat{\mathcal A}
\otimes{\mathcal A}$, the comultiplication on $A_D$ exactly coincides with the one
on $D(A)$, as given in Definition \ref{qdouble}.

\begin{prop}\label{comultiplicationA_D}
For $\phi\in\hat{\mathcal A}$ and $f\in{\mathcal A}$, we have:
\begin{align}
\Delta_D\bigl(\Pi(\phi\otimes f)\bigr)&=\Delta_D\bigl(\pi'(\lambda_{\phi})\pi(L_f)\bigr)
\notag \\
&=\bigl[(\pi'\otimes\pi')(\hat{\Delta}\phi)\bigr]\bigl[(\pi\otimes\pi)(\Delta f)\bigr]
\notag \\
&=(\Pi\otimes\Pi)\left(\sum \phi_{(1)}\otimes f_{(1)}\otimes\phi_{(2)}\otimes f_{(2)}\right).
\notag
\end{align}
\end{prop}

\begin{proof}
Note that by Definition \ref{defnA_D}, we have:
\begin{align}
\Delta_D\bigl(\pi'(\lambda_{\phi})\pi(L_f)\bigr)&={V_D}^*\bigl(1\otimes1\otimes
\pi'(\lambda_{\phi})\pi(L_f)\bigr)V_D  \notag  \\
&=\bigl[{V_D}^*\bigl(1\otimes1\otimes\pi'(\lambda_{\phi})\bigr)V_D\bigr]
\bigl[{V_D}^*\bigl(1\otimes1\otimes\pi(L_f)\bigr)V_D\bigr].  \notag
\end{align}
But by definition of $V_D$ and by definition of $\pi'$, we have:
\begin{align}
{V_D}^*\bigl(1\otimes1\otimes\pi'(\lambda_{\phi})\bigr)V_D&=X_{13}^*Z_{12}Y_{24}^*
Z_{12}^*(1\otimes1\otimes\lambda_{\phi}\otimes1)Z_{12}Y_{24}Z_{12}^*X_{13}  \notag \\
&=X_{13}^*(1\otimes1\otimes\lambda_{\phi}\otimes1)X_{13}  \notag \\
&=\bigl[(\lambda\otimes\lambda)(\hat{\Delta}\phi)\bigr]_{13}
=(\pi'\otimes\pi')(\hat{\Delta}\phi).  \notag
\end{align}
Similarly,
\begin{align}
{V_D}^*\bigl(1\otimes1\otimes\pi(L_f)\bigr)V_D&=X_{13}^*Z_{12}Y_{24}^*Z_{12}^*
\bigl[Z_{34}(1\otimes L_f)_{34}Z_{34}^*\bigr]Z_{12}Y_{24}Z_{12}^*X_{13}  \notag \\
&=Z_{12}Z_{34}Y_{24}^*Z_{12}^*X_{13}^*\bigl[(1\otimes L_f)\bigr]_{34}X_{13}Z_{12}
Y_{24}Z_{34}^*Z_{12}^*  \notag \\
&=Z_{12}Z_{34}Y_{24}^*\bigl[(1\otimes1\otimes1\otimes L_f)\bigr]Y_{24}Z_{34}^*Z_{12}^*
\notag \\
&=Z_{12}Z_{34}\bigl[(L\otimes L)(\Delta f)\bigr]_{24}Z_{34}^*Z_{12}^*
=(\pi\otimes\pi)(\Delta f). \notag
\end{align}
In the second equality above, we used the result of Lemma \ref{lemmaV}\,(3).

Combining these results, we prove the proposition.
\end{proof}

\begin{cor}
The maps $\pi:A\to M(A_D)$ and $\pi':\hat{A}^{\operatorname{op}}\to M(A_D)$, as
defined earlier, are $C^*$-bialgebra homomorphisms.
\end{cor}

\begin{proof}
We already know from Corollary of Proposition \ref{A_D} that $\pi$ and $\pi'$
are $C^*$-algebra homomorphisms.  Meanwhile, from the proof of Proposition
\ref{comultiplicationA_D}, we see that:
$$
(\pi\otimes\pi)\circ\Delta=\Delta_D\circ\pi,\qquad {\text { and}}\qquad
(\pi'\otimes\pi')\circ\hat{\Delta}=\Delta_D\circ\pi'.
$$
\end{proof}

 Propositions \ref{Pi}, \ref{prodA_D}, \ref{comultiplicationA_D} support our
assertion that $(A_D,\Delta_D)$ is indeed the $C^*$-algebraic analog of the
quantum double $D(A)=\hat{A}^{\operatorname{op}}\Join A$, given in Definition
\ref{qdouble}.  To continue with our construction, we next define the antipodal
map $S_D$ on $A_D$.

\begin{lem}\label{ZJhatJ}
With the notation as in Section 2, we have:
$$
Z(J\otimes\hat{J})=(J\otimes\hat{J})Z,\qquad{\text {and}}\qquad
Z(\hat{J}\otimes J)=(\hat{J}\otimes J)Z^*.
$$
\end{lem}

\begin{rem}
The results of the lemma can be shown by a straightforward calculation.
Only the first result is immediately needed, but the second result will be
useful in Section 4.
\end{rem}

\begin{prop}\label{S_D}
Let $\widehat{J_D}$ be the involutive operator in ${\mathcal B}({\mathcal H}
\otimes{\mathcal H})$ defined by $\widehat{J_D}:= J\otimes\hat{J}$.  Then let
$S_D:A_D\to A_D$ by
$$
S_D(x):=\widehat{J_D}x^*\widehat{J_D}=(J\otimes\hat{J})x^*(J\otimes\hat{J}).
$$
In particular, if $f\in{\mathcal A}$ and $\phi\in\hat{\mathcal A}$, we have:
$$
S_D\bigl(\Pi(\phi\otimes f)\bigr)=S_D\bigl(\pi'(\lambda_{\phi})\pi(L_f)\bigr)
=\pi\bigl(S(L_f)\bigr)\pi'\bigl(\hat{S}(\lambda_{\phi})\bigr).
$$
This defines the ``antipode'' on $A_D$.  It is an anti-automorphism on $A_D$,
satisfying: ${S_D\bigl({S_D(x)}^*\bigr)}^*=x$ and $(S_D\otimes S_D)\bigl
(\Delta_D(x)\bigr)={\Delta_D}^{\operatorname{cop}}(S_D(x))$, for $x\in A_D$.
Here ${\Delta_D}^{\operatorname {cop}}=\chi_{1\leftrightarrow3}^{2\leftrightarrow4}
\circ\Delta_D$, where $\chi$ denotes the flip.
\end{prop}

\begin{proof}
For $f\in{\mathcal A}$ and $\phi\in\hat{\mathcal A}$,
\begin{align}
S_D\bigl(\pi'(\lambda_{\phi})\pi(L_f)\bigr)&=(J\otimes\hat{J})\bigl(\pi'(\lambda_{\phi})
\pi(L_f)\bigr)^*(J\otimes\hat{J})  \notag \\
&=(J\otimes\hat{J})\pi({L_f}^*)\pi'({\lambda_{\phi}}^*)(J\otimes\hat{J})  \notag \\
&=(J\otimes\hat{J})Z(1\otimes{L_f}^*)Z^*({\lambda_{\phi}}^*\otimes1)
(J\otimes\hat{J})  \notag \\
&=Z(1\otimes \hat{J}{L_f}^*\hat{J})(J\otimes\hat{J})Z^*(J\otimes\hat{J})
(J{\lambda_{\phi}}^*J\otimes1)  \notag \\
&=Z(1\otimes \hat{J}{L_f}^*\hat{J})Z^*(J{\lambda_{\phi}}^*J\otimes1)
=\pi\bigl(S(L_f)\bigr)\pi'\bigl(\hat{S}(\lambda_{\phi})\bigr).  \notag
\end{align}
In the fourth and fifth equalities, we used the result of Lemma \ref{ZJhatJ}.
In the last equality, we used the definitions of $S$ and $\hat{S}$, given in
terms of $\hat{J}$ and $J$ (See Section 2.).

By definition, it is easy to see that $S_D$ is an anti-automorphism and also that
${S_D\bigl({S_D(x)}^*\bigr)}^*=x$, for $x\in A_D$.  The last statement is also
easy to verify, remembering the corresponding properties of $S$ and $\hat{S}$.
For instance, for $f\in{\mathcal A}$ and $\phi\in\hat{\mathcal A}$, we have:
\begin{align}
&(S_D\otimes S_D)\bigl(\Delta_D(\Pi(\phi\otimes f))\bigr)  \notag \\
&=(S_D\otimes S_D)\bigl([(\pi'\otimes\pi')(\hat{\Delta}\phi)][(\pi\otimes\pi)
(\Delta f)]\bigr)  \notag \\
&=(\pi\otimes\pi)\bigl((S\otimes S)(\Delta f)\bigr)(\pi'\otimes\pi')
\bigl((\hat{S}\otimes\hat{S})(\hat{\Delta}\phi)\bigr)  \notag  \\
&=\bigl[(\pi\otimes\pi)\bigl(\Delta^{\operatorname{cop}}(S(f))\bigr)\bigr]
\bigl[(\pi'\otimes\pi')\bigl(\hat{\Delta}^{\operatorname{cop}}(\hat{S}(\phi))
\bigr)\bigr]  \notag \\
&={\Delta_D}^{\operatorname{cop}}\bigl(\pi(S(f))\pi'(\hat{S}(\phi))\bigr)
={\Delta_D}^{\operatorname{cop}}\bigl(S_D(\Pi(\phi\otimes f))\bigr).
\notag
\end{align}
The third equality follows from the properties of $S$ (in Proposition 2.4 of \cite
{BJKppha}) and of $\hat{S}$ (in Proposition 2.4 of \cite{BJKqhg}).  The fourth
equality follows from Proposition \ref{comultiplicationA_D}.
\end{proof}

 For $S_D$ to be correctly considered an antipode of $(A_D,\Delta_D)$, we further need
the notion of Haar weight clarified.  This will be done later in this paper (in Section 5).
But for our immediate purposes, result of Proposition~\ref{S_D} is sufficient.  In fact,
it is not difficult to show that the definition of $S_D$ given above is equivalent to
the map
$$
S_D:(\operatorname{id}\otimes\operatorname{id}\otimes\Omega)(V_D)\mapsto
(\operatorname{id}\otimes\operatorname{id}\otimes\Omega)({V_D}^*),
$$
for $\Omega\in{\mathcal B}({\mathcal H}\otimes{\mathcal H})_*$.  The general theory
assures us that once we establish the existence of the Haar weight, this map actually
characterizes the antipode (See \cite{BS}, \cite{Wr7}, \cite{KuVa}.).

 According to the general theory of locally compact quantum groups \cite{KuVa},
the antipode allows the ``polar decomposition''.  In our case, with $S_D$ being an
anti-automorphism, its polar decomposition is trivial: That is, $S_D\equiv R_D$
(the ``unitary antipode''), and $\tau_D\equiv\operatorname{Id}$ (the ``scaling
group'').  These observations, in addition to the fact that ${S_D}^2\equiv
\operatorname{Id}$, manifests that $A_D$ is essentially a Kac $C^*$-algebra
(in the sense of \cite{Va}).  This is to be expected, since $(A,\Delta)$ and
$(\hat{A},\hat{\Delta})$ are also as such.

\section{The dual of the quantum double}

 In the previous section, we considered the {\em quantum double\/} $(A_D,\Delta_D)$,
together with its ``antipode'' $S_D$, all within the $C^*$-algebra framework.  The
discussion on its Haar weight (thereby establishing it as a locally compact quantum
group) will be postponed until Section 5.  In the present section, we will consider
the dual object of $(A_D,\Delta_D)$.

 As we can expect from the way $A_D$ was constructed in Definition \ref{defnA_D}
(via the multiplicative unitary operator $V_D$), the dual object will be obtained
by considering the ``left slices'' of $V_D$, as follows (See again, \cite{BS},
\cite{Wr7}.):

\begin{defn}\label{defnhatA_D}
Let $\widehat{A_D}$ be the $C^*$-algebra contained in ${\mathcal B}({\mathcal H}
\otimes{\mathcal H})$, defined by
$$
\widehat{A_D}=\overline{\bigl\{(\Omega\otimes\operatorname{id}\otimes\operatorname
{id})(V_D):\Omega\in{\mathcal B}({\mathcal H}\otimes{\mathcal H})_*\bigr\}}^{\|\ \|}.
$$
In addition, define the {\em comultiplication\/} $\widehat{\Delta_D}:
\widehat{A_D}\to M(\widehat{A_D}\otimes\widehat{A_D})$ by
$$
\widehat{\Delta_D}(y):={V_D}(y\otimes1\otimes1){V_D}^*,\qquad {\text {for $y\in
\widehat{A_D}$}},
$$
which is a non-degenerate $C^*$-homomorphism satisfying the coassociativity
condition: $(\widehat{\Delta_D}\otimes\operatorname{id})\widehat{\Delta_D}
=(\operatorname{id}\otimes\widehat{\Delta_D})\widehat{\Delta_D}$.  Moreover,
$\widehat{\Delta_D}(\widehat{A_D})(\widehat{A_D}\otimes1)$ and $\widehat
{\Delta_D}(\widehat{A_D})(1\otimes\widehat{A_D})$ are dense subsets in
$\widehat{A_D}\otimes\widehat{A_D}$.
\end{defn}

 By the multiplicativity of $V_D$, we know that $(\widehat{A_D},\widehat{\Delta_D})$
is a $C^*$-bialgebra, dual to $(A_D,\Delta_D)$.  Let us now find a more explicit
realization of $(\widehat{A_D},\widehat{\Delta_D})$.  The proof is adapted from
Proposition 8.14 of \cite{BS}. 

\begin{prop}\label{hatA_D}
As a $C^*$-algebra, we have: $\widehat{A_D}=A\otimes\hat{A}^{\operatorname{op}}$.
\end{prop}

\begin{proof}
Given $\Omega\in{\mathcal B}({\mathcal H}\otimes{\mathcal H})_*$ and for
arbitrary $a\in A$, $b\in\hat{A}^{\operatorname{op}}$, define $\tilde
{\Omega}\in{\mathcal B}({\mathcal H}\otimes{\mathcal H})_*$ by $\tilde
{\Omega}:=(1\otimes b)\Omega(a\otimes1)Z^*$.  In particular, if $\Omega
=\Omega_{\xi,\eta}$ (This is a standard notation: $\Omega_{\xi,\eta}(T)
=\langle T\xi,\eta\rangle$, for $T\in{\mathcal B}({\mathcal H}\otimes
{\mathcal H})$ and $\xi,\eta\in{\mathcal H}\otimes{\mathcal H}$.), then
we will have: $\widetilde{\Omega_{\xi,\eta}}(T)=\bigl\langle(a\otimes1)
Z^*T(1\otimes b)\xi,\eta\bigr\rangle$.  With the new notation, we have:
\begin{align}
(\tilde{\Omega}\otimes\operatorname{id}\otimes\operatorname{id})(V_D)
&=(\tilde{\Omega}\otimes\operatorname{id}\otimes\operatorname{id})
(Z_{12}Y_{24}Z_{12}^*X_{13})  \notag \\
&=(\Omega\otimes\operatorname{id}\otimes\operatorname{id})\bigl(
(a\otimes1\otimes1\otimes1)Y_{24}Z_{12}^*X_{13}(1\otimes b\otimes1
\otimes1)\bigr)   \notag \\
&=(\Omega\otimes\operatorname{id}\otimes\operatorname{id})\bigl(
Y_{24}(a\otimes1\otimes1\otimes1)Z_{12}^*(1\otimes b\otimes1\otimes1)
X_{13}\bigr)  \notag \\
&=(\Omega\otimes\operatorname{id}\otimes\operatorname{id})\bigl(
Y_{24}z_{12}X_{13}\bigr),
\notag
\end{align}
where $z=(a\otimes1)Z^*(1\otimes b)$.  Note that since $Z^*=Y\widehat
{\widehat{Y}}^*$, with $Y\in M(A\otimes\hat{A}^{\operatorname{op}})$
and $\widehat{\widehat{Y}}\in M(A^{\operatorname{op}}\otimes\hat{A})$
being elements of (two-sided) multiplier algebras, we see easily that
$z\in A\otimes\hat{A}^{\operatorname{op}}$.  It follows that
$$
\widehat{A_D}=\overline{\bigl\{(\Omega\otimes\operatorname{id}\otimes
\operatorname{id})\bigl(Y_{24}(f\otimes\phi\otimes1\otimes1)X_{13}\bigr):
\Omega\in{\mathcal B}({\mathcal H}\otimes{\mathcal H})_*, f\in A, \phi
\in\hat{A}^{\operatorname{op}}\bigr\}}^{\|\ \|}.
$$

Now note that $\overline{\bigl\{(k\otimes1)Y(\phi\otimes1):\phi\in\hat
{A}^{\operatorname{op}}, k\in{\mathcal K}({\mathcal H})\bigr\}}^{\|\ \|}
={\mathcal K}({\mathcal H})\otimes\hat{A}^{\operatorname{op}}$, and
similarly $\overline{\bigl\{(f\otimes1)X(k\otimes1): f\in A,k\in{\mathcal K}
({\mathcal H})\bigr\}}^{\|\ \|}={\mathcal K}({\mathcal H})\otimes A$.
These follow from the fact that $Y\in M(A\otimes\hat{A}^{\operatorname{op}})
\subseteq M\bigl({\mathcal K}({\mathcal H})\otimes\hat{A}^{\operatorname{op}}
\bigr)$ and $X\in M(\hat{A}^{\operatorname{op}} \otimes A)\subseteq M\bigl
({\mathcal K}({\mathcal H})\otimes A\bigr)$, together with the result
observed in Proposition~3.5\,(2) of \cite{BJKqhg},where we saw that:
$\hat{A}^{\operatorname{op}}=\overline{\bigl\{(\omega\otimes\operatorname
{id})(Y):\omega\in{\mathcal B}({\mathcal H})_*\bigr\}}^{\|\ \|}$ and that
$A=\overline{\bigl\{(\omega\otimes\operatorname{id})(X):\omega\in{\mathcal B}
({\mathcal H})_*\bigr\}}^{\|\ \|}$. [Recall that $Y=\widehat{U_A}$ and
$X=\Sigma{\widehat{U_A}}^*\Sigma$.]

It follows that: $\widehat{A_D}=A\otimes\hat{A}^{\operatorname{op}}$.
\end{proof}

 The proposition shows that in the case of $\widehat{A_D}$, by being just
the (ordinary) tensor product $A\otimes\hat{A}^{\operatorname{op}}$, there
is no ``twisting'' in the algebra structure (Recall that in the case of $A_D$,
it is the coalgebra structure that does not involve twisting.  See Definition
\ref{qdouble} and Proposition \ref{comultiplicationA_D}.).

 On the other hand, we see below that the comultiplication on $\widehat{A_D}$
is equivalent to a ``$\tau$-tensor product'', where $\tau$ is an ``inversion''.
Recall first the definition of an inversion: See Definitions 8.1, 8.2 and
Proposition 8.3 of \cite{BS} (See also \cite{Mj1}.).

\begin{defn}\label{inversion}
\begin{enumerate}
\item Let $(A,\delta_A)$ and $(B,\delta_B)$ be two $C^*$-bialgebras.  An
{\em inversion\/} on $A$ and $B$ is a ${}^*$-isomorphism $\tau:A\otimes B
\to B\otimes A$ such that:
$$
(\tau\otimes\operatorname{id}_A)(\operatorname{id}_A\otimes\tau)(\delta_A
\otimes\operatorname{id}_B)=(\operatorname{id}_B\otimes\delta_A)\tau
$$
and
$$
(\operatorname{id}_B\otimes\tau)(\tau\otimes\operatorname{id}_B)(\operatorname
{id}_A\otimes\delta_B)=(\delta_B\otimes\operatorname{id}_A)\tau,
$$
where we used the same notation $\tau$ for its extension to the multiplier
algebra $M(A\otimes B)$.
\item Given an inversion $\tau$ on $(A,\delta_A)$ and $(B,\delta_B)$, we can
define the map $\delta_{\tau}:A\otimes B\to M(A\otimes B\otimes A\otimes B)$ by
$$
\delta_{\tau}:=(\operatorname{id}_A\otimes\tau\otimes\operatorname{id}_B)
(\delta_A\otimes\delta_B).
$$
Then $\delta_{\tau}$ is coassociative.  It is called the comultiplication
associated with $\tau$.
\end{enumerate}
\end{defn}

 In our case, we can show that the operator $Z$ provides an inversion on
$(A,\Delta^{\operatorname{cop}})$ and $(\hat{A}^{\operatorname{op}},
\hat{\Delta}^{\operatorname{cop}})$, where $\Delta^{\operatorname{cop}}$
and $\hat{\Delta}^{\operatorname{cop}}$ are co-opposite comultiplications.

\begin{prop}\label{tauinversion}
The map $\tau:p\mapsto\Sigma ZpZ^*\Sigma$, where $\Sigma$ is the flip,
is an ``inversion'' on $(A,\Delta^{\operatorname{cop}})$ and
$(\hat{A}^{\operatorname{op}},\hat{\Delta}^{\operatorname{cop}})$.
\end{prop}

\begin{proof}
Let $p\in A\otimes\hat{A}^{\operatorname{op}}$.  Note that $Y^*pY\in 
A\otimes\hat{A}^{\operatorname{op}}$, since $Y=\widehat{U_A}\in
M(A\otimes\hat{A}^{\operatorname{op}})$.  Note also that $p\widehat
{\widehat{Y}}=\widehat{\widehat{Y}}p$, since $\widehat{\widehat{Y}}
\in M(A^{\operatorname{op}}\otimes\hat{A})$.  Therefore,
$$
\tau(p)=\Sigma ZpZ^*\Sigma=\Sigma Y^*\widehat{\widehat{Y}}p{\widehat
{\widehat{Y}}}^*Y\Sigma=\Sigma Y^*pY\Sigma\in\hat{A}^{\operatorname{op}}
\otimes A.
$$
Since $Z$ is a unitary operator, it is clear that $\tau:A\otimes
\hat{A}^{\operatorname{op}}\to\hat{A}^{\operatorname{op}}\otimes A$
is a ${}^*$-isomorphism.

To verify that $\tau$ is an inversion, we note that
$\Delta^{\operatorname{cop}}(a)=X(a\otimes1)X^*$, for $a\in A$,
and that $\hat{\Delta}^{\operatorname{cop}}(b)=Y(b\otimes1) Y^*$,
for $b\in\hat{A}^{\operatorname{op}}$ (These are consequences of
Proposition 3.5\,(2) of \cite{BJKqhg}.).  Indeed, we have:
\begin{align}
&(\operatorname{id}_B\otimes\tau)(\tau\otimes\operatorname{id}_B)
(\operatorname{id}_A\otimes\delta_B)(p)  \notag \\
&=\Sigma_{23}Y_{23}^*\Sigma_{12}Y_{12}^*Y_{23}(p_{12})Y_{23}^*Y_{12}
\Sigma_{12}Y_{23}\Sigma_{23}  \notag \\
&=\Sigma_{23}\Sigma_{12}Y_{13}^*Y_{12}^*Y_{23}(p_{12})Y_{23}^*Y_{12}
Y_{13}\Sigma_{12}\Sigma_{23} \notag \\
&=\Sigma_{23}\Sigma_{12}Y_{23}Y_{12}^*(p_{12})Y_{12}Y_{23}^*
\Sigma_{12}\Sigma_{23} \notag \\
&=\Sigma_{23}Y_{13}\Sigma_{12}Y_{12}^*(p_{12})Y_{12}\Sigma_{12}Y_{13}^*
\Sigma_{23}=\Sigma_{23}Y_{13}\bigl[\tau(p)\bigr]_{12}Y_{13}^*\Sigma_{23}
\notag \\
&=Y_{12}\Sigma_{23}\bigl[\tau(p)\bigr]_{12}\Sigma_{23}Y_{12}^*
=Y_{12}\bigl[\tau(p)\bigr]_{13}Y_{12}^*  \notag \\
&=(\delta_B\otimes\operatorname{id}_A)\tau(p).  \notag
\end{align}
For convenience, we let $B=\hat{A}^{\operatorname{op}}$ and $\delta_B
=\hat{\Delta}^{\operatorname{cop}}$.  The first equality is by applying
the definitions given above, and in the third equality, we used the fact
that $Y$ is multiplicative (i.\,e. $Y_{12}Y_{13}Y_{23}=Y_{23}Y_{12}$ is
equivalent to $Y_{13}^*Y_{12}^*=Y_{23}Y_{12}^*Y_{23}^*$).

A similar computation will verify the other condition, thereby giving us
the proof that $\tau$ is an inversion on $(A,\Delta^{\operatorname{cop}})$
and $(\hat{A}^{\operatorname{op}},\hat{\Delta}^{\operatorname{cop}})$.
\end{proof}

\begin{rem}
This is actually re-writing the proof of Theorem 8.17 in \cite{BS}.
We nevertheless chose to carry out the explicit computation here
(instead of just referring), so that we can have a more tangible
description of the twisted comultiplication below.
\end{rem}

\begin{cor}
On $A\otimes\hat{A}^{\operatorname{op}}$, we have a coassociative
comultiplication $\delta_{\tau}$ defined by
$$
\delta_{\tau}:
=(\operatorname{id}_A\otimes\tau\otimes\operatorname{id}_{\hat{A}^{\operatorname{op}}})
(\Delta^{\operatorname{cop}}\otimes\hat{\Delta}^{\operatorname{cop}}).
$$
Moreover, for $p\in A\otimes\hat{A}^{\operatorname{op}}$, we have:
$\delta_{\tau}(p)=V(p\otimes1\otimes1)V^*$, where $V=Z_{12}^*X_{13}Z_{12}Y_{24}$
is as defined in Lemma \ref{lemmaV}.
\end{cor}

\begin{proof}
The first part is an immediate consequence of Proposition \ref{tauinversion}.
Direct proof is possible for the second part (using similar method as in the
above proof), but we will instead refer the reader to the proof of Lemma 8.9
in \cite{BS}.
\end{proof}

 We are now ready to give a more specific description of the (twisted)
comultiplication on $\widehat{A_D}$.

\begin{prop}\label{comultiplicationhatA_D}
For $y\in\widehat{A_D}$, we have:
$$
\widehat{\Delta_D}(y)=(Z\otimes Z)\bigl[\delta_{\tau}(Z^*yZ)\bigr](Z^*\otimes Z^*).
$$
\end{prop}

\begin{proof}
Write $y=ZpZ^*$, for $p\in A\otimes\hat{A}^{\operatorname{op}}$. [Note that for
any $y\in\widehat{A_D}=A\otimes\hat{A}^{\operatorname {op}}$, we have: $Z^*yZ
\in A\otimes\hat{A}^{\operatorname{op}}$.] Then:
\begin{align}
\widehat{\Delta_D}(y)&=V_D(y\otimes1\otimes1){V_D}^*
=V_D(ZpZ^*\otimes1\otimes1){V_D}^*  \notag \\
&=(Z\otimes Z)V(Z^*\otimes Z^*)[ZpZ^*\otimes1\otimes1](Z\otimes Z)V^*(Z^*\otimes Z^*)
\notag \\
&=(Z\otimes Z)V(p\otimes1\otimes1)V^*(Z^*\otimes Z^*)  \notag \\
&=(Z\otimes Z)\bigl[\delta_{\tau}(p)\bigr](Z^*\otimes Z^*)
=(Z\otimes Z)\bigl[\delta_{\tau}(Z^*yZ)\bigr](Z^*\otimes Z^*).  \notag
\end{align}
The third equality uses Lemma \ref{lemmaV} (See also the proof of the Corollary
following it, where it is noted that $V_D=(Z\otimes Z)V(Z^*\otimes Z^*)$.).  The
next to last equality follows from the Corollary of Proposition \ref{tauinversion}
above.
\end{proof}

\begin{prop}\label{hatA_Disom}
We have the following $C^*$-bialgebra isomorphisms:
\begin{align}
(\widehat{A_D},\widehat{\Delta_D})&\cong(A\otimes\hat{A}^{\operatorname{op}},
\delta_{\tau}),\qquad {\text {where $\delta_{\tau}=(\operatorname{id}\otimes
\tau\otimes\operatorname{id})(\Delta^{\operatorname{cop}}\otimes
\hat{\Delta}^{\operatorname{cop}})$,}}   \notag \\
&\cong(A^{\operatorname{op}}\otimes\hat{A},\delta_{\tau'}),
\qquad {\text {where $\delta_{\tau'}=(\operatorname{id}
\otimes{\tau}'\otimes\operatorname{id})(\Delta\otimes\hat{\Delta})$.}}   
\notag
\end{align}
Here $\tau:A\otimes\hat{A}^{\operatorname{op}}\to\hat{A}^{\operatorname{op}}
\otimes A$ is as above, and ${\tau}':A^{\operatorname{op}}\otimes\hat{A}
\to\hat{A}\otimes A^{\operatorname{op}}$ is defined by ${\tau}'(q)=\Sigma 
Z^*qZ\Sigma$, which is also an inversion.
\end{prop}

\begin{proof}
The first isomorphism is the result of Proposition \ref{comultiplicationhatA_D},
given by the map $\widehat{A_D}\ni y\mapsto Z^*yZ\in A\otimes\hat{A}^{\operatorname
{op}}$.  The second isomorphism is given by the map $A\otimes\hat{A}^{\operatorname
{op}}\ni p\mapsto (S\otimes\hat{S})(p)\in A^{\operatorname{op}}\otimes\hat{A}$
(Note that in our case, the antipodes $S$ and $\hat{S}$ are both anti-automorphisms.).

To verify the last description of $\delta_{\tau'}$, let $q\in A^{\operatorname{op}}
\otimes\hat{A}$.  Then:
\begin{align}
\delta_{\tau'}(q)&=(S\otimes\hat{S}\otimes S\otimes\hat{S})
\bigl[\delta_{\tau}\bigl((S\otimes\hat{S})(q)\bigr)\bigr] \notag \\
&=(S\otimes\hat{S}\otimes S\otimes\hat{S})\bigl[(\operatorname{id}
\otimes\tau\otimes\operatorname{id})(\Delta^{\operatorname{cop}}\otimes
\hat{\Delta}^{\operatorname{cop}})\bigl((S\otimes\hat{S})(q)\bigr)\bigr]
\notag \\
&=(S\otimes\hat{S}\otimes S\otimes\hat{S})\bigl[(\operatorname{id}
\otimes\tau\otimes\operatorname{id})(S\otimes S\otimes\hat{S}\otimes\hat{S})
\bigl((\Delta\otimes\hat{\Delta})(q)\bigr)\bigr]   \notag \\
&=\bigl(\operatorname{id}\otimes(\hat{S}\otimes S)\tau(S\otimes\hat{S})
\otimes\operatorname{id}\bigr)\bigl((\Delta\otimes\hat{\Delta})(q)\bigr).
\notag
\end{align}
But for any $q\in A^{\operatorname{op}}\otimes\hat{A}$, we have:
\begin{align}
(\hat{S}\otimes S)\tau(S\otimes\hat{S})(q)&=(\hat{S}\otimes S)\tau\bigl(
(\hat{J}\otimes J)q^*(\hat{J}\otimes J)\bigr)  \notag \\
&=(\hat{S}\otimes S)\bigl[\Sigma Z\bigl((\hat{J}\otimes J)q^*(\hat{J}\otimes J)
\bigr)Z^*\Sigma\bigr]  \notag \\
&=(J\otimes\hat{J})\bigl[\Sigma Z\bigl((\hat{J}\otimes J)q^*(\hat{J}\otimes J)
\bigr)Z^*\Sigma\bigr]^*(J\otimes\hat{J})  \notag \\
&=\Sigma(\hat{J}\otimes J)Z(\hat{J}\otimes J)q(\hat{J}\otimes J)Z^*(\hat{J}
\otimes J)\Sigma  \notag \\
&=\Sigma Z^*qZ\Sigma={\tau}'(q).   \notag
\end{align}
The next to last equality uses Lemma \ref{ZJhatJ}.  By this,
we have proven the second isomorphism.
\end{proof}

\begin{rem}
Note that we may also regard $\tau'$ as an inversion on $(A,\Delta)$
and $(\hat{A}^{\operatorname{op}},\hat{\Delta})$ [The definition given
in the above proposition is still valid.].  Then we will have: $\tau'
=\chi{\tau}^{-1}\chi$, where $\chi$ is the flip.  In this setting, the
comultiplication $\delta_{\tau'}$ is equivalent to the co-opposite
comultiplication $\delta_{\tau}^{\operatorname{cop}}$, in the sense
that for any $p\in A\otimes\hat{A}^{\operatorname{op}}$, we have:
\begin{equation}\label{(deltataucop)}
(Z^*\otimes Z^*)\bigl[\delta_{\tau'}(ZpZ^*)\bigr](Z\otimes Z)
=\delta_{\tau}^{\operatorname{cop}}(p)
=\chi_{1\leftrightarrow3}^{2\leftrightarrow4}\bigl(\delta_{\tau}(p)\bigr).
\end{equation}
This may be shown by direct computation (which we do not carry out here),
but it is really a consequence of our working with the Kac systems \cite{BS}.
Then by definition of $\widehat{\Delta_D}$, it follows from equation \eqref
{(deltataucop)} that:
$$
\widehat{\Delta_D}^{\operatorname{cop}}(y)
=(Z\otimes Z)\bigl[\delta_{\tau}^{\operatorname{cop}}(Z^*yZ)\bigr](Z^*\otimes Z^*)
=\delta_{\tau'}(y),
$$
giving us a very tidy description of $\widehat{\Delta_D}^{\operatorname{cop}}
\equiv\delta_{\tau'}$.
\end{rem}

 Propositions \ref{comultiplicationhatA_D} and \ref{hatA_Disom} provide us
a specific description of the comultiplication on $\widehat{A_D}$, and we
see that $\widehat{\Delta_D}$ is equivalent to a certain ``twisted'' tensor
product comultiplication.  Now, let us look for the antipodal map on
$(\widehat{A_D},\widehat{\Delta_D})$.  As the general theory suggests and
similar to the cases of $S$, $\hat{S}$, and $S_D$ (see \cite{BS}, \cite{Wr7},
\cite{KuVa}, as well as \cite{BJKppha}, \cite{BJKqhg}, and the paragraph
following Proposition \ref{S_D}), we wish to consider the following map:
\begin{equation}\label{(hatS_D)}
\widehat{S_D}:(\Omega\otimes\operatorname{id}\otimes\operatorname{id})
(V_D)\mapsto(\Omega\otimes\operatorname{id}\otimes\operatorname{id})
({V_D}^*),\qquad\Omega\in{\mathcal B}({\mathcal H}\otimes{\mathcal H})_*.
\end{equation}
But note that by Lemma \ref{lemmaV}\,(3),
\begin{align}
{V_D}^*&=(Z_{12}Y_{24}Z_{12}^*X_{13})^*=(Z_{12}Z_{34}Z_{12}^*X_{13}Z_{12}Y_{24}Z_{34}^*
Z_{12}^*)^*  \notag \\
&=Z_{12}Z_{34}Y_{24}^*Z_{12}^*X_{13}^*Z_{34}^*=Z_{34}(Z_{12}Y_{24}^*Z_{12}^*X_{13}^*)
Z_{34}^*.  \notag
\end{align}
From this, the equation \eqref{(hatS_D)} can be written as:
$$
\widehat{S_D}:(\Omega\otimes\operatorname{id}\otimes\operatorname{id})
(Z_{12}Y_{24}Z_{12}^*X_{13})\mapsto Z\bigl[(\Omega\otimes\operatorname{id}
\otimes\operatorname{id})(Z_{12}Y_{24}^*Z_{12}^*X_{13}^*)\bigr]Z^*.
$$
Remembering the characterizations of $S:(\omega\otimes\operatorname{id})(X)
\mapsto(\omega\otimes\operatorname{id})(X^*)$ and $\hat{S}:(\omega'\otimes
\operatorname{id})(Y)\mapsto(\omega'\otimes\operatorname{id})(Y^*)$, the
new expression suggests that
$$
\widehat{S_D}:a\otimes b\mapsto Z\bigl(S(a)\otimes\hat{S}(b)\bigr)Z^*,\qquad
a\in A,b\in\hat{A}^{\operatorname{op}}.
$$
Having this result as a motivation, we treat this more precisely in the
following proposition:

\begin{prop}\label{hatS_D}
Let $J_D$ be defined by $J_D:=Z(\hat{J}\otimes J)\in{\mathcal B}({\mathcal H}
\otimes{\mathcal H})$.  By Lemma \ref{ZJhatJ}, we have: $J_D=Z(\hat{J}\otimes J)
=(\hat{J}\otimes J)Z^*$, and it follows that $J_D$ is involutive.  We then
define $\widehat{S_D}:\widehat{A_D}\to\widehat{A_D}$ by
$$
\widehat{S_D}(y):=J_D y^* J_D=Z(\hat{J}\otimes J)y^*(\hat{J}\otimes J)Z^*
=(\hat{J}\otimes J)Z^*y^*Z(\hat{J}\otimes J).
$$
Then $\widehat{S_D}$ defines the ``antipode'' on $\widehat{A_D}$.
\end{prop}

\begin{proof}
Note first that the definition given in the proposition is same as the one
suggested in the previous paragraph:
$$
\widehat{S_D}(y)=J_D y^*J_D
=Z(\hat{J}\otimes J)y^*(\hat{J}\otimes J)Z^*=Z\bigl((S\otimes\hat{S})(y)
\bigr)Z^*.
$$
In the last equality, we used the definitions of the maps $S$ and $\hat{S}$
as given in Section 2.

From its definition, it is immediate that $\widehat{S_D}$ is an anti-automorphism
on $\widehat{A_D}$, satisfying: $\widehat{S_D}\bigl(\widehat{S_D}(y)^*\bigr)^*
=y$, for $y\in\widehat{A_D}$.  Next, let us also prove:
\begin{equation}\label{(antipodehatS_D)}
(\widehat{S_D}\otimes\widehat{S_D})\bigl(\widehat{\Delta_D}(y)\bigr)
={\widehat{\Delta_D}}^{\operatorname{cop}}\bigl(\widehat{S_D}(y)\bigr)
=\chi_{1\leftrightarrow3}^{2\leftrightarrow4}\bigl[\widehat{\Delta_D}
(\widehat{S_D}(y))\bigr].
\end{equation}
For this, let us consider without loss of generality $y=ZpZ^*=Z(a\otimes b)Z^*$,
where, $a\in A$ and $b\in\hat{A}^{\operatorname{op}}$.  Then by the remark
following Proposition~\ref{hatA_Disom}, the right hand side of equation
\eqref{(antipodehatS_D)} can be realized as being equal to $\delta_{\tau'}
\bigl(\widehat{S_D}(y)\bigr)=\delta_{\tau'}\bigl((S\otimes\hat{S}) (p)\bigr)$.
Meanwhile, the left hand side can be written as follows:
\begin{align}
&(\widehat{S_D}\otimes\widehat{S_D})\bigl(\widehat{\Delta_D}(y)\bigr)  \notag \\
&=(\hat{J}\otimes J\otimes\hat{J}\otimes J)(Z^*\otimes Z^*)\bigl[\widehat{\Delta_D}
(y)^*\bigr](Z\otimes Z)(\hat{J}\otimes J\otimes\hat{J}\otimes J)  \notag \\
&=(\hat{J}\otimes J\otimes\hat{J}\otimes J)\bigl[\delta_{\tau}(p)^*\bigr]
(\hat{J}\otimes J\otimes\hat{J}\otimes J)  \notag \\
&=(S\otimes\hat{S}\otimes S\otimes\hat{S})\bigl[\delta_{\tau}(a\otimes b)\bigr] 
\notag \\
&=(S\otimes\hat{S}\otimes S\otimes\hat{S})\bigl[(\operatorname{id}\otimes
\tau\otimes\operatorname{id})(\Delta^{\operatorname{cop}}a\otimes
\hat{\Delta}^{\operatorname{cop}}b)\bigr]  \notag \\
&=(\operatorname{id}\otimes\tau'\otimes\operatorname{id})\bigl[(S\otimes S
\otimes\hat{S}\otimes\hat{S})\bigl(\Delta^{\operatorname{cop}}a\otimes
\hat{\Delta}^{\operatorname{cop}}b\bigr)\bigr]  \notag \\
&=(\operatorname{id}\otimes\tau'\otimes\operatorname{id})\bigl[\Delta\bigl
(S(a)\bigr)\otimes\hat{\Delta}\bigl(\hat{S}(b)\bigr)\bigr]  \notag \\
&=\delta_{\tau'}\bigl(S(a)\otimes\hat{S}(b)\bigr)=\delta_{\tau'}\bigl((S\otimes
\hat{S})(p)\bigr).
\notag
\end{align}
In this way, we verify the equation \eqref{(antipodehatS_D)}.  Note that in the
second equality above, we used the characterization of $\widehat{\Delta_D}$ given
in Proposition \ref{comultiplicationhatA_D}.  In the fifth equality, we used the
result: $(\hat{S}\otimes S)\bigl[\tau(p)\bigr]=\tau'\bigl[(S\otimes\hat{S})(p)
\bigr]$, which appeared in the proof of Proposition \ref{hatA_Disom}.  The sixth
equality is using the properties of the antipode maps $S$ and $\hat{S}$. 
\end{proof}

 As in the case of $(A_D,\Delta_D)$ and its antipode $S_D$, the properties of
the antipodal map $\widehat{S_D}$ noted above (including $\widehat{S_D}^2\equiv
\operatorname{Id}$) manifests that $(\widehat{A_D},\widehat{\Delta_D})$ is
again a Kac $C^*$-algebra (with the existence of its Haar weight, to be given
below).

\section{Haar weight}

To show that the quantum double $(A_D,\Delta_D)$ and its dual $(\widehat{A_D},
\widehat{\Delta_D})$ constructed above are indeed locally compact quantum groups,
we need a discussion on their Haar weights.  These Haar weights will be described
in terms of the Haar weights on $(A,\Delta)$ and $(\hat{A},\hat{\Delta})$, which
we recall below:

\begin{lem}\label{varphi}
On ${\mathcal A}$, define a linear functional $\varphi$ by
$$
\varphi(a)=\int a(0,0,r)\,dr.
$$
It can be extended to a faithful, lower semi-continuous, tracial weight (still
denoted by $\varphi$) on the $C^*$-algebra $A$. The weight $\varphi$ satisfies
the ``left invariance property'': For any positive element $a\in A$ such that
$\varphi(a)<\infty$, and for $\omega\in A^*_+$, we have:
$$
\varphi\bigl((\omega\otimes\operatorname{id})(\Delta a)\bigr)=\omega(1)
\varphi(a).
$$
\end{lem}

\begin{lem}\label{hatvarphi}
On $\hat{\mathcal A}$, define a linear functional $\hat{\varphi}$ by
$$
\hat{\varphi}(b)=\int b(x,y,0)\,dxdy.
$$
It can be extended to a faithful, lower semi-continuous, KMS weight (still
denoted by $\hat{\varphi}$) on the $C^*$-algebra $\hat{A}$.  The weight
$\hat{\varphi}$ also satisfies the ``left invariance property'': For any
positive element $b\in\hat{A}$ such that $\hat{\varphi}(b)<\infty$, and
for positive $\omega\in{\hat{A}}^*$, we have:
$$
\hat{\varphi}\bigl((\omega\otimes\operatorname{id})(\hat{\Delta}b)\bigr)
=\omega(1)\hat{\varphi}(b).
$$
Moreover, we have the following (unimodular) property: $\hat{\varphi}\circ
\hat{S}=\hat{\varphi}$.
\end{lem}

\begin{rem}
In both cases above, $\omega(1)=\|\omega\|$.  These results are described
in \cite{BJKppha} (section~3) and in \cite{BJKqhg} (section~2), respectively.
In particular, the left invariance properties are proved in Theorem 3.9 of
\cite{BJKppha} and in Theorem 2.11 of \cite{BJKqhg}. The invariance properties
written above are in a weak form, but the general theory assures us that
they are actually sufficient (See \cite{KuVa}, and see also section~1 of
\cite{BJKppha}.).  Note finally that $\hat{\varphi}$ is invariant under
$\hat{S}$ (so $(\hat{A},\hat{\Delta})$ is unimodular), while it is not the
case for $\varphi$ and $S$.
\end{rem}

 The invariance properties stay the same when we consider instead
$(A^{\operatorname{op}},\Delta)$ and $(\hat{A}^{\operatorname {op}},\hat
{\Delta})$.  On the other hand, see the Corollary below for the cases having
co-opposite comultiplications.

\begin{cor}
Since the antipode $S:A\to A$ satisfies $(S\otimes S)\Delta=\Delta^{\operatorname
{cop}}\circ S$, it follows from Lemma \ref{varphi} that $\varphi\circ S$ is left
invariant for $\Delta^{\operatorname{cop}}$.  That is,
$$
(\varphi\circ S)\bigl((\omega\otimes\operatorname{id})(\Delta^{\operatorname{cop}}a)
\bigr)=\omega(1)\varphi\bigl(S(a)\bigr).
$$
Similarly, from Lemma \ref{hatvarphi} and by using the property of the antipode
$\hat{S}$, we have:
$$
\hat{\varphi}\bigl((\omega\otimes\operatorname{id})(\hat{\Delta}^{\operatorname
{cop}}b)\bigr)=(\hat{\varphi}\circ\hat{S})\bigl((\omega\otimes\operatorname{id})
(\hat{\Delta}^{\operatorname{cop}}b)\bigr)=\omega(1)\hat{\varphi}\bigl(S(b)\bigr)
=\omega(1)\hat{\varphi}(b).
$$
\end{cor}

\begin{proof}
The verification is straightforward.  In the second case, we are using the fact
that $\hat{\varphi}$ is invariant under $\hat{S}$.
\end{proof}

 Let us begin our discussion by considering the Haar weight on $(\widehat{A_D},
\widehat{\Delta_D})$, whose definition is given below. [There is actually a simpler
characterization for $\widehat{\varphi_D}$, as can be found in Corollary of Proposition
\ref{hatvarphiunimodular} below.  But our choice of the definition has been made so
that the proofs of the later propositions are a little simpler.]
 
\begin{defn}\label{Haarweightdual}
Let $\widehat{\varphi_D}$ be the faithful, lower semi-continuous weight on
$\widehat{A_D}$, defined by
$$
\widehat{\varphi_D}(y):=\bigl((\varphi\circ S)\otimes\hat{\varphi}\bigr)(Z^*yZ).
$$
\end{defn}

\begin{rem}
Since $\varphi\circ S$ and $\hat{\varphi}$ are densely defined weights on
the $C^*$-algebras $A$ and $\hat{A}^{\operatorname{op}}$, and since they are
both faithful and lower semi-continuous (i.\,e. they are ``proper'' weights),
we can define their tensor product $(\varphi\circ S)\otimes\hat{\varphi}$.
See Definition~1.27 of \cite{KuVa}. Therefore, $\widehat{\varphi_D}$ is a proper
weight on the $C^*$-algebra $\widehat{A_D}$.  See the following lemma.
\end{rem}

\begin{lem}\label{tproductweight}
For arbitrary proper weights $\varphi$ and $\psi$ on $C^*$-algebras $A$ and
$B$, consider the tensor product weight $\varphi\otimes\psi$ on $A\otimes B$.
Then the set ${\mathfrak N}_{\varphi}\odot{\mathfrak N}_{\psi}(\subseteq
{\mathfrak N}_{\varphi\otimes\psi})$ forms a core for the GNS map
$\Lambda_{\varphi\otimes\psi}$.
\end{lem}

\begin{rem}
For a more systematic discussion on tensor product weights, see section 1.6
and Appendix of \cite{KuVa}.  The result of this lemma actually goes back
to Haagerup's density theorem \cite{H}: Adapted to our case, if $X\in
{\mathfrak N}_{\varphi\otimes\psi}$, there exists a sequence $\{X_n\}$
in ${\mathfrak N}_{\varphi}\odot{\mathfrak N}_{\psi}$ such that:
$\lim_{n\to\infty}\Lambda_{\varphi\otimes\psi}(X_n)=\Lambda_{\varphi
\otimes\psi}(X)$.  Because of this density result, pretty much all
the properties of a tensor product weight $\varphi\otimes\psi$ can be
obtained by just working with the elements from the algebraic tensor
product ${\mathfrak N}_{\varphi}\odot{\mathfrak N}_{\psi}$.  
\end{rem}

 The weight $\widehat{\varphi_D}$ is shown to be left invariant for
$(\widehat{A_D},\widehat{\Delta_D})$, which is actually unimodular.  See
Proposition \ref{hatvarphiunimodular} and Theorem \ref{hatleftinvariance}
below.

\begin{prop}\label{hatvarphiunimodular}
The weight $\widehat{\varphi_D}$ is $\widehat{S_D}$-invariant: $\widehat
{\varphi_D}=\widehat{\varphi_D}\circ\widehat{S_D}$.
\end{prop}

\begin{proof}
Let $y\in\widehat{A_D}$ be such that $y=Z(L_f\otimes\lambda_{\phi})Z^*$,
where $f\in{\mathcal A}$ and $\phi\in\hat{\mathcal A}$.  Since such elements
are dense in $\widehat{A_D}$ (and form a core for $\widehat{\varphi_D}$),
our proof will be achieved if we just verify that $\widehat{\varphi_D}(y)
=(\widehat{\varphi_D}\circ\widehat{S_D})(y)$.

By definition, we can see easily that:
\begin{align}
\widehat{\varphi_D}(y)&=\bigl((\varphi\circ S)\otimes\hat{\varphi}\bigr)(Z^*yZ)
=(\varphi\circ S)(L_f)\hat{\varphi}(\lambda_{\phi})  \notag \\
&=\left(\int(e^{2\lambda r})^n f(0,0,-r)\,dr\right)\left(\int\phi(x,y,0)\,dxdy
\right)  \notag \\
&=\int(e^{-2\lambda r})^n f(0,0,r)\phi(x,y,0)\,dxdydr.
\notag
\end{align}

Meanwhile, by Proposition \ref{hatS_D}, we know that $\widehat{S_D}(y)=S(L_f)
\hat{S}(\lambda_{\phi})$.  Therefore,
$$
\widehat{\varphi_D}\bigl(\widehat{S_D}(y)\bigr)=\bigl((\varphi\circ S)\otimes
\hat{\varphi}\bigr)\bigl(Z^*S(L_f)\hat{S}(\lambda_{\phi})Z\bigr).
$$
To compute this (so that we can compare the result with $\widehat{\varphi_D}(y)$
given above), let us begin by finding a suitable realization of $Z^*S(L_f)\hat{S}
(\lambda_{\phi})Z$.  Let $\xi\in{\mathcal H}\otimes{\mathcal H}$ and compute:
\begin{align}
&Z^*S(L_f)\hat{S}(\lambda_{\phi})Z\xi(x,y,r;x',y',r')  \notag \\
&=\int(e^{-\lambda r})^n\bar{e}\bigl[\eta_{\lambda}(r)\beta(x-e^{-\lambda r}x',
e^{-\lambda r}y')\bigr]e\bigl[\eta_{\lambda}(r)\beta(e^{\lambda r'-\lambda r}x',
y-e^{-\lambda r}y')\bigr]  \notag \\
&\qquad S(f)(\tilde{x},\tilde{y},r)\hat{S}(\phi)(e^{\lambda\tilde{r}-\lambda r}x',
e^{\lambda\tilde{r}-\lambda r}y',r'-\tilde{r})  \notag \\
&\qquad\bar{e}\bigl[\eta_{\lambda}(r)\beta(\tilde{x},
y+e^{\lambda r'-\lambda r}y'-e^{-\lambda r}y'-\tilde{y})\bigr]  \notag \\
&\qquad\bar{e}\bigl[\eta_{\lambda}(r)\beta(e^{\lambda\tilde{r}-\lambda r}x',
y+e^{\lambda r'-\lambda r}y'-e^{-\lambda r}y'-\tilde{y}-e^{\lambda\tilde{r}
-\lambda r}y')\bigr]  \notag \\
&\qquad(e^{\lambda r})^n e\bigl[\eta_{\lambda}(r)\beta(x+e^{\lambda r'-\lambda r}x'
-e^{-\lambda r}x'-\tilde{x}-e^{\lambda\tilde{r} -\lambda r}x',e^{-\lambda r}y')\bigr]
\notag \\
&\qquad\xi(x+e^{\lambda r'-\lambda r}x'-\tilde{x}-e^{\lambda\tilde{r}-\lambda r}x',
y+e^{\lambda r'-\lambda r}y'-\tilde{y}-e^{\lambda\tilde{r}-\lambda r}y',r;x',y',
\tilde{r})\,d\tilde{x}d\tilde{y}d\tilde{r}  \notag \\
&=\int F(\tilde{x},\tilde{y},r;e^{\lambda\tilde{r}}x',e^{\lambda\tilde{r}}y',
\tilde{r}-r')  \notag \\
&\qquad\bar{e}\bigl[\eta_{\lambda}(r)\beta(\tilde{x},y-\tilde{y})\bigr]
\xi(x-\tilde{x},y-\tilde{y},r;x',y',\tilde{r})\,d\tilde{x}d\tilde{y}d\tilde{r}
\notag \\
&=(L\otimes\lambda)_F\xi(x,y,r;x',y',r'),
\notag
\end{align}
where $F$ is defined by
\begin{align}
F(\tilde{x},\tilde{y},r;x',y',\tilde{r})&=(e^{2\lambda r})^n f(-e^{\lambda r}
\tilde{x}-e^{-\lambda\tilde{r}}x'+x',-e^{\lambda r}\tilde{y}-e^{-\lambda\tilde{r}}y'+y',
-r)  \notag \\
&\quad\bar{e}\bigl[\eta_{\lambda}(r)\beta(\tilde{x}+e^{-\lambda\tilde{r}-\lambda r}x'
-e^{-\lambda r}x',e^{-\lambda\tilde{r}-\lambda r}y')\bigr]
\notag \\
&\quad\bar{e}\bigl[\eta_{\lambda}(-\tilde{r})\beta(e^{-\lambda r}x',e^{-\lambda r}y')
\bigr]e\bigl[\eta_{\lambda}(r)\beta(e^{-\lambda r}x',\tilde{y})\bigr]
\notag \\
&\quad\bar{e}\bigl[\eta_{\lambda}(r)\beta(\tilde{x},\tilde{y})\bigr]
\phi(-e^{-\lambda\tilde{r}-\lambda r}x',-e^{-\lambda\tilde{r}-\lambda r}y',\tilde{r}).
\notag
\end{align}
The expression for $F$ is obtained by remembering the definitions of $S(f)$ and
$\hat{S}(\phi)$ given in Section 2, and by using the change of variables.  In this
way, we found out the following realization:
$$
Z^*S(L_f)\hat{S}(\lambda_{\phi})Z=(L\otimes\lambda)_F.
$$
This means that:
\begin{align}
\widehat{\varphi_D}\bigl(\widehat{S_D}(y)\bigr)&=\bigl((\varphi\circ S)\otimes
\hat{\varphi}\bigr)(F) \notag \\
&=\int (e^{-2\lambda r})^n F(0,0,r;x',y',0)\,dx'dy'dr  \notag \\
&=\int f(0,0,-r)\phi(-e^{-\lambda r}x',-e^{-\lambda r}y',0)\,dx'dy'dr  \notag \\
&=\int f(0,0,r)\phi(x',y',0)(e^{-2\lambda r})^n\,dx'dy'dr.
\notag
\end{align}

Combining the results, we see that $\widehat{\varphi_D}\bigl(\widehat{S_D}(y)\bigr)
=\widehat{\varphi_D}(y)$, proving our assertion.
\end{proof}

\begin{cor}
For any $a\in A$ and $b\in\hat{A}^{\operatorname{op}}$, we have:
$$
\bigl((\varphi\circ S)\otimes\hat{\varphi}\bigr)\bigl(Z^*(a\otimes b)Z\bigr)
=\varphi(a)\hat{\varphi}\bigl(\hat{S}(b)\bigr)=\varphi(a)\hat{\varphi}(b).
$$
In particular, for any $p\in\widehat{A_D}$, we have: $\widehat{\varphi_D}(p)
=(\varphi\otimes\hat{\varphi})(p)$. [This is the simpler characterization of
$\widehat{\varphi_D}$ mentioned earlier.]
\end{cor}

\begin{proof}
Consider $y=Z\bigl(S(a)\otimes\hat{S}(b)\bigr)Z^*$.  Then by Proposition
\ref{hatS_D}, we know that $\widehat{S_D}(y)=S\bigl(S(a)\bigr)\otimes\hat{S}
\bigl(\hat{S}(b)\bigr)=a\otimes b$.  It follows that:
$$
\widehat{\varphi_D}\bigl(\widehat{S_D}(y)\bigr)=\bigl((\varphi\circ S)\otimes
\hat{\varphi}\bigr)(Z^*(a\otimes b)Z\bigr).
$$
On the other hand, by definition of $\widehat{\varphi_D}$ and by using the
unimodularity of $\hat{\varphi}$, we have:
$$
\widehat{\varphi_D}(y)=\bigl((\varphi\circ S)\otimes\hat{\varphi}\bigr)
\bigl(S(a)\otimes\hat{S}(b)\bigr)=\varphi(a)\hat{\varphi}\bigl(\hat{S}(b)\bigr)
=\varphi(a)\hat{\varphi}(b).
$$
Since we should have $\widehat{\varphi_D}\bigl(\widehat{S_D}(y)\bigr)=\widehat
{\varphi_D}(y)$ by Proposition \ref{hatvarphiunimodular}, the first statement
follows.  The second statement is an immediate consequence.
\end{proof}

\begin{theorem}\label{hatleftinvariance}
For any positive element $y\in\widehat{A_D}$ such that $\widehat{\varphi_D}(y)
<\infty$, and for $\Omega\in\widehat{A_D}^*_+$, we have:
$$
\widehat{\varphi_D}\bigl((\Omega\otimes\operatorname{id}\otimes
\operatorname{id})(\widehat{\Delta_D}(y))\bigr)=\Omega(1)\widehat
{\varphi_D}(y).
$$ 
\end{theorem}

\begin{proof}
For convenience, let us write $B=\hat{A}^{\operatorname{op}}$.  Consider
$y=Z(a\otimes b)Z^*$, where $a\in A_+$, $(\varphi\circ S)(a)<\infty$ and
$b\in B_+$, $\hat{\varphi}(b)<\infty$. Assume also that $\Omega$ has the
form $\Omega=\omega_1\otimes\omega_2$, for some $\omega_1\in A^*_+$ and
$\omega_2\in B^*_+$.  Then:
\begin{align}
&\bigl((\Omega\otimes\operatorname{id}\otimes\operatorname{id})(\widehat{\Delta_D}(y))
\bigr)   \notag \\
&=(\omega_1\otimes\omega_2\otimes\operatorname{id}\otimes\operatorname{id})\bigl(
(Z_{12}Z_{34})\bigl[\Sigma_{23}Z_{23}(\Delta^{\operatorname{cop}}a\otimes\hat
{\Delta}^{\operatorname{cop}}b)Z^*_{23}\Sigma_{23}\bigr](Z_{12}Z_{34})^*\bigr)
\notag \\
&=(\tilde{\omega_1}\otimes\tilde{\omega_2}\otimes\operatorname{id}\otimes
\operatorname{id})\bigl(Z_{34}\Sigma_{23}Z_{23}(\Delta^{\operatorname{cop}}a
\otimes\hat{\Delta}^{\operatorname{cop}}b)Z^*_{23}\Sigma_{23}Z^*_{34}\bigr),
\notag
\end{align}
where $\tilde{\omega_1}$, $\tilde{\omega_2}$ are defined such that $(\tilde{\omega_1}
\otimes\tilde{\omega_2})(\cdot):=(\omega_1\otimes\omega_2)(Z\cdot Z^*)$.

But by Lemma \ref{ZY^*}, and by using $\hat{\Delta}^{\operatorname{cop}}b
=Y(b\otimes1)Y^*$, we have:
\begin{align}
&Z_{34}\Sigma_{23}Z_{23}(\Delta^{\operatorname{cop}}a\otimes\hat
{\Delta}^{\operatorname{cop}}b)Z^*_{23}\Sigma_{23}Z^*_{34}  \notag \\
&=Y^*_{34}\Sigma_{23}Y^*_{23}\bigl[(\operatorname{id}\otimes\operatorname{id}
\otimes\hat{\Delta}^{\operatorname{cop}})(\Delta^{\operatorname{cop}}a\otimes b)
\bigr]Y_{23}\Sigma_{23}Y_{34}  \notag \\
&=\Sigma_{23}Y^*_{24}Y^*_{23}Y_{34}(\Delta^{\operatorname{cop}}a\otimes b
\otimes1)Y^*_{34}Y_{23}Y_{24}\Sigma_{23}  \notag \\
&=\Sigma_{23}Y_{34}Y^*_{23}(\Delta^{\operatorname{cop}}a\otimes b\otimes1)
Y_{23}Y^*_{34}\Sigma_{23}  \notag \\
&=\Sigma_{23}\bigl[(\operatorname{id}\otimes\operatorname{id}\otimes\hat
{\Delta}^{\operatorname{cop}})(Y^*_{23}(\Delta^{\operatorname{cop}}a\otimes b)
Y_{23})\bigr]\Sigma_{23}.
\notag
\end{align}
We are using the multiplicativity of $Y$ (i.\,e. $Y_{23}Y_{24}Y_{34}=Y_{34}Y_{23}$)
in the third equality.  Therefore, by using the definition of $\widehat{\varphi_D}$
given in Definition \ref{Haarweightdual} and in Corollary of Proposition \ref
{hatvarphiunimodular}, we have:
\begin{align}
&\widehat{\varphi_D}\bigl((\Omega\otimes\operatorname{id}\otimes\operatorname{id})
(\widehat{\Delta_D}(y))\bigr)=(\varphi\otimes\hat{\varphi})\bigl((\Omega\otimes
\operatorname{id}\otimes\operatorname{id})(\widehat{\Delta_D}(y))\bigr)  \notag \\
&=(\varphi\otimes\hat{\varphi})\bigl((\tilde{\omega_1}\otimes\operatorname{id}
\otimes\tilde{\omega_2}\otimes\operatorname{id})\bigl[(\operatorname{id}
\otimes\operatorname{id}\otimes\hat{\Delta}^{\operatorname{cop}})(Y^*_{23}
(\Delta^{\operatorname{cop}}a\otimes b)Y_{23})\bigr]\bigr).
\notag
\end{align}

Further computation shows the following:
\begin{align}
&\widehat{\varphi_D}\bigl((\Omega\otimes\operatorname{id}\otimes\operatorname{id})
(\widehat{\Delta_D}(y))\bigr)  \notag \\
&=\tilde{\omega_2}(1)(\varphi\otimes\hat{\varphi})\bigl((\tilde{\omega_1}\otimes
\operatorname{id})\bigl[Y^*_{23}(\Delta^{\operatorname{cop}}a\otimes b)Y_{23}\bigr]
\bigr)  \notag \\
&=\tilde{\omega_2}(1)\bigl((\varphi\circ S)\otimes\hat{\varphi}\bigr)\bigl((\tilde
{\omega_1}\otimes\operatorname{id})[\Delta^{\operatorname{cop}}a\otimes b]\bigr)
\notag \\
&=\tilde{\omega_2}(1)\tilde{\omega_1}(1)(\varphi\circ S)(a)\hat{\varphi}(b)
\notag \\
&=\omega_2(1)\omega_1(1)\widehat{\varphi_D}(y)=\Omega(1)\widehat{\varphi_D}(y).
\notag
\end{align}
Notice that the first and third equalities above use the left invariance properties
of $\hat{\varphi}$ and $(\varphi\circ S)$ (See Corollary of Lemmas \ref{varphi}
and \ref{hatvarphi}.).  However, some care has to be given (using Lemma \ref
{tproductweight}), if we want them to be perfectly valid.  We will not go into the
details here (to avoid our discussion from becoming too technical and lengthy),
but for instance, we may follow the discussion similar to the proof of Lemma~3.5
and Corollary~3.6 of \cite{Ya}.  Meanwhile, the second equality is due to Lemma
\ref{ZY^*} and Corollary of Proposition \ref{hatvarphiunimodular}.  Fourth and
fifth equalities follow from the observation that
\begin{align}
\tilde{\omega_1}(1)\tilde{\omega_2}(1)&=(\tilde{\omega_1}\otimes\tilde{\omega_2})
(1\otimes1)=(\omega_1\otimes\omega_2)\bigl(Z(1\otimes1)Z^*\bigr)  \notag \\
&=(\omega_1\otimes\omega_2)(1\otimes1)=\Omega(1).  \notag
\end{align}

So far we proved the case when $y=Z(a\otimes b)Z^*$ and $\Omega=\omega_1\otimes
\omega_2$.  Extending the proof for general $y\in\widehat{A_D}_+$ and $\Omega
\in\widehat{A_D}^*_+$ is not necessarily trivial.  Nevertheless, we will again
invoke Lemma~\ref{tproductweight} here and refer the reader instead to the papers
mentioned above (See \cite{Ya}, \cite{KuVa}.).
\end{proof}

 Theorem \ref{hatleftinvariance} establishes the proof that $\widehat{\varphi_D}$
is a legitimate (invariant) Haar weight.  By general theory \cite{KuVa}, it is
therefore the unique (up to multiplication by a scalar) Haar weight for $(\widehat
{A_D},\widehat{\Delta_D})$.  In our case, we note that even if $(A,\Delta)$ was
non-unimodular, $\widehat{\varphi_D}$ is actually unimodular for $(\widehat{A_D},
\widehat{\Delta_D})$ (Proposition \ref{hatvarphiunimodular}).  Since this is
the case, we do not need any further discussion on the ``modular function''.
Summarizing the results so far, we now state the following theorem:

\begin{theorem}
The $C^*$-bialgebra $(\widehat{A_D},\widehat{\Delta_D})$, together with its
additional structure maps including the antipode $\widehat{S_D}$ and the
(unimodular) Haar weight $\widehat{\varphi_D}$, is a $C^*$-algebraic locally
compact quantum group, in the sense of Kustermans and Vaes.
\end{theorem}

 As we have made our case throughout Section 4 and Section 5, we regard
$(\widehat{A_D},\widehat{\Delta_D})$ as the dual of the quantum double.
Namely, $\widehat{D(A)}$.

 As for the quantum double $D(A)=(A_D,\Delta_D)$, remembering that it is
the dual object of $(\widehat{A_D},\widehat{\Delta_D})$ associated with
the multiplicative unitary operator $V_D$, and that $(\widehat{A_D},\widehat
{\Delta_D})$ is a legitimate locally compact quantum group (Theorem \ref
{hatleftinvariance}), we conclude immediately from general theory \cite{KuVa},
\cite{Wr7} that it is also a $C^*$-algebraic locally compact quantum group.
This achieves our stated goal.

 For the remainder of this section, let us just give an explicit description
of the Haar weight $\varphi_D$ of $(A_D,\Delta_D)$, whose existence (and
uniqueness up to multiplication by a scalar) is assured from the above
observation.  The subalgebra $\hat{\mathcal A}\otimes{\mathcal A}\subseteq
A_D$ forms a core for the Haar weight.

\begin{prop}
For $\Pi(b\otimes a)=\pi'(\lambda_b)\pi(L_a)\in A_D$, where $b\in\hat
{\mathcal A}$ and $a\in{\mathcal A}$, define:
$$
\varphi_D\bigl(\Pi(b\otimes a)\bigr):=\hat{\varphi}(\lambda_b)\varphi(L_a)
=\int (b\otimes a)(x,y,0;0,0,r')\,dxdydr'.
$$
This defines a linear functional on $\hat{\mathcal A}\otimes{\mathcal A}$.
Then we have:
$$
\varphi_D\bigl(\Pi(\phi\otimes f)^*\Pi(b\otimes a)\bigr)
=\hat\varphi({\lambda_{\phi}}^*\lambda_b)\varphi({L_f}^*L_a),
$$
for $b,\phi\in\hat{\mathcal A}$ and for $a,f\in{\mathcal A}$.  Furthermore,
we have:
$$
(\Omega\otimes{\varphi_D})\bigl({\Delta_D}(\Pi(b\otimes a))\bigr)
=\Omega(1){\varphi_D}\bigl(\Pi(b\otimes a)\bigr),\qquad\Omega
\in{\mathcal B}({\mathcal H}\otimes{\mathcal H})_*.
$$

This will characterize the Haar weight on $(A_D,\Delta_D)$.  In other words,
the functional $\varphi_D$ extends to a (unique) $C^*$-algebra weight on
$A_D$, which is left invariant.
\end{prop}

\begin{proof}
Note that
$$
\Pi(\phi\otimes f)^*\Pi(b\otimes a)=\Pi\bigl((\phi\otimes f)^*\times(b\otimes a)),
$$
where the involution and multiplication on $\hat{\mathcal A}\otimes{\mathcal A}$
are as given in equations \eqref{(doubleinvol)} and \eqref{(doubleprod)}.  By
a straightforward computation using the equations \eqref{(doubleinvol)} and
\eqref{(doubleprod)}, we have
\begin{align}
&\int\bigl((\phi\otimes f)^*\times(b\otimes a)\bigr)(x,y,0;0,0,r')\,dxdydr'  \notag \\
&=(\cdots)  \notag \\
&=\int\overline{\phi(x,y,\tilde{r})}b(x,y,\tilde{r})\overline{f(\tilde{x},\tilde{y},r')}
a(\tilde{x},\tilde{y},r')\,d\tilde{x}d\tilde{y}d\tilde{r}dxdydr'  \notag \\
&=\hat{\varphi}(b\times_{\hat{A}}{\phi}^*)\varphi(f^*\times_A a)
=\hat{\varphi}({\lambda_{\phi}}^*\lambda_b)\varphi({L_f}^*L_a).  \notag
\end{align}
It thus follows that: $\varphi_D\bigl(\Pi(\phi\otimes f)^*\Pi(b\otimes a)\bigr)
=\hat\varphi({\lambda_{\phi}}^*\lambda_b)\varphi({L_f}^*L_a)$.

 Using this, we can give $\hat{\mathcal A}\otimes{\mathcal A}$ a left Hilbert
algebra structure.  Moreover, we can show without difficulty that the GNS
Hilbert space for $\varphi_D$ is ${\mathcal H}\otimes{\mathcal H}$, while the
GNS representation is $\Pi$.  Following a standard procedure (see \cite{Cm2},
and also \cite{BJKppha}, \cite{BJKqhg}), we can define a $C^*$-algebra weight
on $A_D$ extending the functional $\varphi_D$ (The extended weight will be
still denoted by $\varphi_D$.).

 Meanwhile, at least at the level of the (dense) subalgebra $\hat{\mathcal A}
\otimes{\mathcal A}$, the verification of the left invariance of $\varphi_D$
is not very difficult.  Note that by Proposition \ref{comultiplicationA_D},
we can write:
\begin{align}
(\Omega\otimes\operatorname{id}\otimes\operatorname{id})\bigl(\Delta_D(\Pi
(b\otimes a))\bigr)&=\sum(\Omega\otimes\operatorname{id}\otimes\operatorname
{id})\bigl((\Pi\otimes\Pi)(b_{(1)}\otimes a_{(1)}\otimes b_{(2)}\otimes a_{(2)})
\bigr)   \notag \\
&=\sum\bigl[\Omega\bigl(\pi'(b_{(1)})\pi(a_{(1)})\bigr)\bigl(\pi'(b_{(2)})
\pi(a_{(2)})\bigr)\bigr],
\notag
\end{align}
where we are using Sweedler's notation for $\hat{\Delta}b$ and $\Delta a$.
And, for convenience, we regard $b=\lambda_b$ and $a=L_a$.  Then:
\begin{align}
(\Omega\otimes{\varphi_D})\bigl({\Delta_D}(\Pi(b\otimes a))\bigr)
&=\sum\bigl[\Omega\bigl(\pi'(b_{(1)})\pi(a_{(1)})\bigr)
{\varphi_D}\bigl(\pi'(b_{(2)})\pi(a_{(2)})\bigr)\bigr]  \notag \\
&=\sum\bigl[\Omega\bigl((b_{(1)}\otimes1)Z(1\otimes a_{(1)})Z^*\bigr)
\hat{\varphi}(b_{(2)}){\varphi}(a_{(2)})\bigr].   \notag
\end{align}
Without loss of generality, assume that $\Omega=\Omega_{\xi,\eta}$, for
$\xi,\eta\in{\mathcal H}\otimes{\mathcal H}$ (following the standard
notation, as appeared in the proof of Proposition \ref{hatA_D}).  Then
the expression becomes:
\begin{align}
&(\Omega\otimes{\varphi_D})\bigl({\Delta_D}(\Pi(b\otimes a))\bigr)
=\sum\bigl[\bigl\langle(b_{(1)}\otimes1)Z(1\otimes a_{(1)})Z^*\xi,
\eta\bigr\rangle\hat{\varphi}(b_{(2)}){\varphi}(a_{(2)})\bigr]
\notag \\
&=\int\sum\bigl[(b_{(1)}\otimes1)Z(1\otimes a_{(1)})Z^*\xi(x,y,r;x',y',r')
\overline{\eta(x,y,r;x',y',r')}  \notag \\
&\qquad\quad b_{(2)}(\tilde{x},\tilde{y},0)a_{(2)}(0,0,\tilde{r})\bigr]\,dxdydr
dx'dy'dr'd\tilde{x}d\tilde{y}d\tilde{r}.
\notag
\end{align}
We can compute this using the formulas we obtained in Section 2 for
$b_{(1)}=\lambda(b_{(1)})$ and $a_{(1)}=L(a_{(1)})$, as well as the
operator $Z$ (obtained in Section 3).  Next, note that $\hat{\Delta}b
=\sum[b_{(1)}\otimes b_{(2)}]$ and that $\Delta a=\sum[a_{(1)}
\otimes a_{(2)}]$, where we can use the equation \eqref{(hatAcomult)}
for $\hat{\Delta b}$ and the equation \eqref{(Acomult)} for $\Delta a$.
Then, by using change of variables, the expression becomes:
\begin{align}
&(\Omega\otimes{\varphi_D})\bigl({\Delta_D}(\Pi(b\otimes a))\bigr)
\notag \\
&=\int b(e^{\lambda r}x+\tilde{x},e^{\lambda r}y+\tilde{y},0)
a(0,0,r'+\tilde{r})  \notag \\
&\qquad\xi(x,y,r;x',y',r')\overline{\eta(x,y,r;x',y',r')}\,dxdydr
dx'dy'dr'd\tilde{x}d\tilde{y}d\tilde{r}  \notag \\
&=\int b(\tilde{x},\tilde{y},0)a(0,0,\tilde{r})\xi(x,y,r;x',y',r')
\overline{\eta(x,y,r;x',y',r')}\,dxdydrdx'dy'dr'd\tilde{x}d\tilde{y}
d\tilde{r}  \notag \\
&=\hat{\varphi}(b)\varphi(a)\Omega_{\xi,\eta}(1)
=\Omega(1)\varphi_D\bigl(\Pi(b\otimes a)\bigr).
\notag
\end{align}

Since we already know the existence of the unique Haar weight from the
discussion preceding the proposition, this invariance property at the dense
subalgebra level is enough to assure us that $\varphi_D$ is indeed the
legitimate Haar weight for $(A_D,\Delta_D)$.
\end{proof}

 Summarizing the results from Section 3 and the discussion on the Haar weight
given here, we conclude the following:

\begin{theorem}
The $C^*$-bialgebra $(A_D,\Delta_D)$, together with the Haar weight ${\varphi_D}$,
is a $C^*$-algebraic locally compact quantum group.  It is the ($C^*$-algebraic)
``quantum double'': $D(A)=\hat{A}^{\operatorname{op}}\Join A$.
\end{theorem}

\begin{rem}
Unlike in the case of $(\widehat{A_D},\widehat{\Delta_D})$ and its Haar weight
$\widehat{\varphi_D}$, we can show easily that $\varphi_D$ is non-unimodular:
That is, $\varphi_D\circ S_D\ne\varphi_D$.  The same modular function operator
for $(A,\Delta)$ (see section 5 of \cite{BJKppha}) will work as the modular
function for $(A_D,\Delta_D)$.
\end{rem}

\section{Quantum Universal $R$-matrix}

 Just as in the case of the purely algebraic framework, our quantum double
$(A_D,\Delta_D)$ is also equipped with a (quasi-triangular) ``quantum universal
$R$-matrix'' type operator.  The definition of a quantum $R$-matrix in the
$C^*$-algebra framework is essentially same as in the more usual, Hopf algebra
or QUE algebra setting (See \cite{Dr}, \cite{CP}, \cite{Mj}, for the usual
definition; And see section 6 of \cite{BJKp2}, for the definition in the
$C^*$-algebra setting.).

 In this section, we will give a brief construction of the operator ${\mathcal R}
\in M(A_D\otimes A_D)$, which will be considered as the ``quantum universal
$R$-matrix'' for the quantum group $(A_D,\Delta_D)$.  Let us begin with
a lemma, which actually follows from Lemma \ref{lemmaV}.  The proof is
adapted from section~8 of \cite{BS}.

\begin{lem}\label{lemmaR}
Let the notations be as before, and let $X$, $Y$, and $Z$ be the operators
defined earlier.  Then we have:
\begin{enumerate}
\item $Z_{12}^*X_{14}Z_{12}Y_{25}Y_{45}=Y_{45}Y_{25}Z_{12}^*X_{14}Z_{12}$
\item $Z_{34}X_{14}Z_{34}^*X_{15}X_{35}=X_{35}X_{15}Z_{34}X_{14}Z_{34}^*$
\end{enumerate}
\end{lem}

\begin{proof}
From Lemma \ref{lemmaV} (3), we know that:
$$
Z_{34}Z_{12}^*X_{13}Z_{12}Y_{24}=Y_{24}Z_{12}^*X_{13}Z_{12}Z_{34}.
$$
It follows that: $Z_{12}^*X_{13}Z_{12}Y_{24}Z_{34}^*=Z_{34}^*Y_{24}Z_{12}^*
X_{13}Z_{12}$.  Remembering that $X\in M(\hat{A}^{\operatorname{op}}
\otimes A)$, $Y\in M(A\otimes\hat{A}^{\operatorname{op}})$, and that
$Z=\widehat{\widehat{Y}}Y^*$, where $\widehat{\widehat{Y}}\in
M(A^{\operatorname{op}}\otimes\hat{A})$, this becomes:
$Z_{12}^*X_{13}Z_{12}Y_{24}Y_{34}=Y_{34}Y_{24}Z_{12}^*X_{13}Z_{12}$.
(The point here is that $\widehat{\widehat{Y}}_{34}^*$ commutes with
all the operators in the equation.)  Certainly, this is equivalent to (1):
$$
Z_{12}^*X_{14}Z_{12}Y_{25}Y_{45}=Y_{45}Y_{25}Z_{12}^*X_{14}Z_{12}.
$$

For (2), recall first that $X=\Sigma Y^*\Sigma$.  Then (1) can be
re-written as:
$$
Z_{12}^*X_{14}Z_{12}X_{52}^*X_{54}^*=X_{54}^*X_{52}^*Z_{12}^*X_{14}Z_{12}.
$$
From this, we have: $X_{52}X_{54}Z_{12}^*X_{14}Z_{12}=Z_{12}^*X_{14}Z_{12}
X_{54}X_{52}$.  So we have: $Z_{12}X_{52}X_{54}Z_{12}^*X_{14}=X_{14}Z_{12}
X_{54}X_{52}Z_{12}^*$, which is same as:
$$
Z_{12}X_{52}Z_{12}^*X_{54}X_{14}=X_{14}X_{54}Z_{12}X_{52}Z_{12}^*.
$$
But this is actually equivalent to (2) [Legs 1,2,4,5 are now considered
as legs 3,4,5,1.]: $Z_{34}X_{14}Z_{34}^*X_{15}X_{35}=X_{35}X_{15}Z_{34}
X_{14}Z_{34}^*$.
\end{proof}

 We are now ready to give the description of our ``quantum $R$-matrix''
operator ${\mathcal R}$.  Again, the definition is a slight modification
of the one considered in section 8 of \cite{BS}. 

\begin{prop}
Let ${\mathcal R}=Z_{34}X_{14}Z_{34}^*$.  The following properties hold.
\begin{enumerate}
\item ${\mathcal R}\in M(A_D\otimes A_D)$.
\item We have: $(\Delta_D\otimes\operatorname{id})({\mathcal R})
={\mathcal R}_{13}{\mathcal R}_{23}$ and $(\operatorname{id}\otimes
\Delta_D)({\mathcal R})={\mathcal R}_{13}{\mathcal R}_{12}$.
\item For any $x\in A_D$, we have: ${\mathcal R}(\Delta_D(x)){\mathcal R}^*
=\Delta_D^{\operatorname{cop}}(x)$.
\item The operator ${\mathcal R}$ satisfies the ``quantum Yang-Baxter equation'':
\linebreak
${\mathcal R}_{12}{\mathcal R}_{13}{\mathcal R}_{23}
={\mathcal R}_{23}{\mathcal R}_{13}{\mathcal R}_{12}$.
\end{enumerate}
\end{prop}

\begin{rem}
In (2) and (4) above, we are viewing ${\mathcal R}$ as an operator contained
in ${\mathcal B}\bigl(({\mathcal H}\otimes{\mathcal H})\otimes({\mathcal H}
\otimes{\mathcal H})\bigr)$.
\end{rem}

\begin{proof}
(1) Recall that $X\in M(\hat{A}^{\operatorname{op}}\otimes A)$.  So by naturally
extending the $C^*$-algebra homomorphisms $\pi'$ and $\pi$ defined in Proposition
\ref{A_D} and its Corollary, we can see that:
$$
{\mathcal R}=Z_{34}X_{14}Z_{34}^*=(\pi'\otimes\pi)(X)\in M(A_D\otimes A_D).
$$

(2) By using the characterization of ${\mathcal R}$ given above, we have:
\begin{align}
(\Delta_D\otimes\operatorname{id})({\mathcal R})&=(\Delta_D\otimes\operatorname
{id})\bigl((\pi'\otimes\pi)(X)\bigr)=(\pi'\otimes\pi'\otimes\pi)\bigl(\hat{\Delta}
\otimes\operatorname{id})(X)\bigr)  \notag \\
&=(\pi'\otimes\pi'\otimes\pi)(X_{12}^*X_{23}X_{12})=(\pi'\otimes\pi'\otimes\pi)
(X_{13}X_{23})  \notag \\
&=\bigl[(\pi'\otimes\pi'\otimes\pi)(X_{13})\bigr]\bigl[(\pi'\otimes\pi'\otimes\pi)
(X_{23})\bigr]={\mathcal R}_{13}{\mathcal R}_{23}.  \notag
\end{align}
The second equality is due to $\Delta_D\circ\pi'=(\pi'\otimes\pi')\circ\hat
{\Delta}$, which was observed in Corollary of Proposition \ref{comultiplicationA_D}.
Third equality is using the fact that $\hat{\Delta}b=X^*(1\otimes b)X$, for $b\in
\hat{A}^{\operatorname{op}}$, while the next equality is the multiplicativity
of $X$.  The next to the last equality is using that $\pi'$ and $\pi$ are
${}^*$-homomorphisms.

Meanwhile, by remembering that $\Delta a=Y^*(1\otimes a)Y$, for $a\in A$, and that
$Y=\Sigma X^*\Sigma$, a similar computation will give us the other equation:
$(\operatorname{id}\otimes\Delta_D)({\mathcal R})={\mathcal R}_{13}{\mathcal R}_{12}$.

(3) Recall that $b=(\operatorname{id}\otimes\omega)(X)\in\hat{A}^{\operatorname{op}}$,
and $a=(\operatorname{id}\otimes\omega')(Y)\in A$, for $\omega,\omega'\in{\mathcal B}
({\mathcal H})_*$, and that these operators generate $\hat{A}^{\operatorname{op}}$
and $A$, respectively. [This result follows from Proposition 3.5\,(2) of \cite{BJKqhg}
and section 6 of \cite{BJKppha}, and was also noted in the proof of Proposition \ref
{A_D}.]

So consider $b=(\operatorname{id}\otimes\omega)(X)$ and compute.  Then:
\begin{align}
{\mathcal R}\bigl[\Delta_D\bigl(\pi'(b)\bigr)\bigr]&={\mathcal R}\bigl[(\pi'\otimes\pi')
(\hat{\Delta} b)\bigr]=(Z_{34}X_{14}Z_{34}^*)\bigl[X_{13}^*(1\otimes1\otimes b\otimes1)
X_{13}\bigr]  \notag \\
&=(\operatorname{id}\otimes\operatorname{id}\otimes\operatorname{id}\otimes\operatorname{id}
\otimes\omega)(Z_{34}X_{14}Z_{34}^*X_{13}^*X_{35}X_{13}) \notag \\
&=(\operatorname{id}\otimes\operatorname{id}\otimes\operatorname{id}\otimes\operatorname{id}
\otimes\omega)(Z_{34}X_{14}Z_{34}^*X_{15}X_{35}) \notag \\
&=(\operatorname{id}\otimes\operatorname{id}\otimes\operatorname{id}\otimes\operatorname{id}
\otimes\omega)(X_{35}X_{15}Z_{34}X_{14}Z_{34}^*) \notag \\
&=(\operatorname{id}\otimes\operatorname{id}\otimes\operatorname{id}\otimes\operatorname{id}
\otimes\omega)(X_{31}^*X_{15}X_{31}Z_{34}X_{14}Z_{34}^*) \notag \\
&=Y_{13}(b\otimes1\otimes1\otimes1)Y_{13}^*Z_{34}X_{14}Z_{34}^*  \notag \\
&=\bigl[(\pi'\otimes\pi')(\hat{\Delta}^{\operatorname{cop}}b)\bigr]{\mathcal R}
=\bigl[{\Delta_D}^{\operatorname{cop}}\bigl(\pi'(b)\bigr)\bigr]{\mathcal R}.
\notag
\end{align}
The first equality is again using $\Delta_D\circ\pi'=(\pi'\otimes\pi')\circ\hat{\Delta}$.
The fourth and sixth equalities follow from the multiplicativity of $X$, while the fifth
equality is by Lemma \ref{lemmaR}\,(2).  The seventh equality is just using $Y=\Sigma X^*
\Sigma$.

Next, consider $a=(\operatorname{id}\otimes\omega')(Y)$ and compute. Then:
\begin{align}
{\mathcal R}\bigl[\Delta_D\bigl(\pi(a)\bigr)\bigr]&={\mathcal R}\bigl[(\pi\otimes\pi)
(\Delta a)\bigr]=(Z_{34}X_{14}Z_{34}^*)\bigl[Z_{34}Z_{12}Y_{24}^*(1\otimes1\otimes 1
\otimes a)Y_{24}Z_{12}^*Z_{34}^*\bigr]  \notag \\
&=(\operatorname{id}\otimes\operatorname{id}\otimes\operatorname{id}\otimes\operatorname{id}
\otimes\omega)(Z_{34}X_{14}Z_{12}Y_{24}^*Y_{45}Y_{24}Z_{12}^*Z_{34}^*) \notag \\
&=(\operatorname{id}\otimes\operatorname{id}\otimes\operatorname{id}\otimes\operatorname{id}
\otimes\omega)(Z_{34}X_{14}Z_{12}Y_{25}Y_{45}Z_{12}^*Z_{34}^*) \notag \\
&=(\operatorname{id}\otimes\operatorname{id}\otimes\operatorname{id}\otimes\operatorname{id}
\otimes\omega)(Z_{34}Z_{12}Y_{45}Y_{25}Z_{12}^*X_{14}Z_{12}Z_{12}^*Z_{34}^*) \notag \\
&=(\operatorname{id}\otimes\operatorname{id}\otimes\operatorname{id}\otimes\operatorname{id}
\otimes\omega)(Z_{34}Z_{12}Y_{42}^*Y_{25}Y_{42}Z_{12}^*X_{14}Z_{34}^*) \notag \\
&=Z_{34}Z_{12}X_{24}(1\otimes a\otimes1\otimes1)X_{24}^*Z_{12}^*Z_{34}^*Z_{34}X_{14}Z_{34}^* 
\notag \\
&=\bigl[(\pi\otimes\pi)(\Delta^{\operatorname{cop}}a)\bigr]{\mathcal R}
=\bigl[{\Delta_D}^{\operatorname{cop}}\bigl(\pi(a)\bigr)\bigr]{\mathcal R}.
\notag
\end{align}
This is done in exactly same way as in the previous case. In particular, the first
equality is using $\Delta_D\circ\pi=(\pi\otimes\pi)\circ\Delta$ (See Corollary of
Proposition~\ref{comultiplicationA_D}), while the fifth equality uses Lemma \ref
{lemmaR}\,(1).

Since $A_D$ is known to be generated by the operators $\pi'(b)\pi(a)$, we conclude
from the previous two results that we have: ${\mathcal R}\bigl[\Delta_D(x)\bigr]
{\mathcal R}^*={\Delta_D}^{\operatorname{cop}}(x)$, for any $x\in A_D$. [Note that
by definition, ${\mathcal R}$ is unitary.]

(4) The last statement is an immediate consequence of results (2) and (3):
$$
{\mathcal R}_{12}{\mathcal R}_{13}{\mathcal R}_{23}={\mathcal R}_{12}
\bigl[(\Delta_D\otimes\operatorname{id})({\mathcal R})\bigr]=\bigl[
({\Delta_D}^{\operatorname{cop}}\otimes\operatorname{id})({\mathcal R})
\bigr]{\mathcal R}_{12}={\mathcal R}_{23}{\mathcal R}_{13}{\mathcal R}_{12}.
$$
First equality follows from (2); the second equality is from (3); and the third
equality is from (2) with the legs 1 and 2 interchanged.
\end{proof}

 Existence of a quantum $R$-matrix for a Hopf algebra (or a quantum group)
is quite useful in the development of the representation theory (See, for
instance, \cite{BJKhj}.).  However, we will postpone to a future occasion
any further discussion about the operator ${\mathcal R}$ and its applications.
Some of these future discussions will be about the relationship between $(A_D,
\Delta_D)$ and its ``classical limit'', which is the double Poisson--Lie
group considered in \cite{BJKhj}.


\bibliography{refqd}

\bibliographystyle{amsplain}

\end{document}